\documentclass[A4j,11pt]{article}
\usepackage{amsfonts,amsmath,amscd}
\usepackage{graphics}
\usepackage{amssymb}
\begin{document}

\title{{\bf The complexifications of pseudo-Riemannian\\
manifolds and anti-Kaehler geometry}}
\author{{\bf Naoyuki Koike}}

\date{}

\maketitle

%
%
%

%
\begin{abstract}
In this paper, we first define 
the complexification of a real analytic map between real analytic Koszul 
manifolds and show that the complexified map is the holomorphic extension of 
the original map.  Next we define an anti-Kaehler metric compatible with 
the adapted complex structure on the complexification of a real analytic 
pseudo-Riemannian manifold.  
In particular, for a pseudo-Riemannian homogeneous space, we define another 
complexification and a (complete) anti-Kaehler metric on 
the complexification.  
One of main purposes of this paper is to find the interesting relation between 
these two complexifications (equipped with the anti-Kaehler metrics) of a 
pseudo-Riemannian homogeneous space.  
Another of main purposes of this paper is to show that almost all principal 
orbits of some isometric action on the first complexification (equipped with 
the anti-Kaehler metric) of a semi-simple pseudo-Riemannian symmetric space 
are curvature-adapted isoparametric submanifolds with flat section in the 
sense of this paper.  
\end{abstract}




\section{Introduction}
Any $C^{\omega}$-manifold $M$ admit its complexification, that is, 
a complex manifold equipped with an anti-holomorphic involution $\sigma$ whose 
fixed point set is $C^{\omega}$-diffeomorphic to $M$, where 
$C^{\omega}$ means the real analyticity.  To get a canonical complexification of $M$ one needs 
some extra structure on $M$.  For example, if $M$ equips with a $C^{\omega}$-Riemannian metric $g$, 
then so-called adapted complex structure $J^g$ is 
defined on a tubular neighborhood $U^g$ (which we take as largely as possible) 
of the zero section of the tangent bundle $TM$ of $M$ and 
$(U^g,J^g)$ gives a 
complexification of $M$ under the identification of $M$ with the zero section 
(see [21,25]).  We denote $(U^g,J^g)$ by $M_g^{\Bbb C}$.  In more general, 
R. Sz$\ddot o$ke ([30]) extended the 

\vspace{0.5truecm}


\noindent
$\overline{\qquad\qquad\qquad\qquad\qquad\qquad}$

\vspace{0.1truecm}

\noindent
{\footnotesize 2010 {\sl Mathematics Subject classification.} Primary 53C56; Secondly 53C42.

\noindent
{\sl Key words and phrases.} complexification, adapted complex structure, 
anti-Kaehler metric,

\noindent
isoparametric submanifold.}

\newpage

\noindent
notion of the adapted complex structure 
to the case where $M$ equips with a $C^{\omega}$-Koszul connection $\nabla$, 
where a $C^{\omega}$-Koszul connection means 
a $C^{\omega}$-linear connection of $TM$.  In this paper, we denote 
this complex structure by $J^{\nabla}$, its domain by $U^{\nabla}$ and 
$(U^{\nabla},J^{\nabla})$ by $M_{\nabla}^{\Bbb C}$, which is a complexification 
of $M$.  We shall call a manifold equipped with a Koszul 
connection a {\it Koszul manifold}.  Thus we get a canonical complexification 
of a $C^{\omega}$-Koszul manifold (as a special case, a 
$C^{\omega}$-pseudo-Riemannian manifold).  On the complexification 
$M_g^{\Bbb C}:=(U^g,J^g)$ of a $C^{\omega}$-pseudo-Riemannian manifold $(M,g)$ 
of index $\nu$, a pseudo-Kaehler metric $g_K$ of index $\nu$ compatible with 
$J^g$ which satisfies $\iota^{\ast}g_K=\frac12 g$ ($\iota\,:\,$ the inclusion 
map of $M$ into $M_g^{\Bbb C}$) is defined in terms of the energy function 
$E:TM\to{\Bbb R}$ (see [30] in detail), where $E$ is defined by 
$E(v):=\frac12g(v,v)\,\,(v\in TM)$.  

In [15], we defined the (extrinsic) complexification of a complete 
$C^{\omega}$-Riemannian submanifold $(M,g)$ immersed by $f$ in a Riemannian 
symmetric space $N=G/K$ of non-compact type as follows.  
First we defined the 
complexification $f^{\Bbb C}$ of $f$ as a map of a tubular neighborhood 
$(M_g^{\Bbb C})_f$ of $M$ in the complexification $M_g^{\Bbb C}$ of $M$ into 
the anti-Kaehler symmetric space $G^{\Bbb C}/K^{\Bbb C}$.  Next we showed that 
$f^{\Bbb C}$ is an immersion over a tubular neighborhood 
$(M_g^{\Bbb C})_{f,i}$ of the zero section in $M_g^{\Bbb C}$.  
We called an anti-Kaehler submanifold 
$((M_g^{\Bbb C})_{f,i},(f^{\Bbb C}\vert_{(M_g^{\Bbb C})_{f,i}})^{\ast}
\langle\,\,,\,\,\rangle)$ in $G^{\Bbb C}/K^{\Bbb C}$ the {\it extrinsic 
complexification} of the Riemannian submanifold $(M,g)$.  Also, in [15], 
we showed that complex focal radii of $M$ introduced in [14] are the 
quantities which indicate the position of focal points of 
$((M_g^{\Bbb C})_{f,i},(f^{\Bbb C}\vert_{(M_g^{\Bbb C})_{f,i}})^{\ast}
\langle\,\,,\,\,\rangle)$.  
Furthermore, by imposing a condition related to complex focal radii, 
we defined the notions of a complex equifocal submanifold and proper complex 
equifocal submanifold.  
It is conjectured that this notion coincides with that of an isoparametric 
submanifold with flat section introduced by Heintze-Liu-Olmos in [11].  
In [15], [16] and [17], we obtained some results for a complex equifocal 
submanifold by investigating the lift of the complexification of 
the submanifold to some path space.  

L. Geatti and C. Gorodski [6] showed that a polar representation 
of a real reductive algebraic group on a pseudo-Euclidean space has the same 
closed orbits as the isotropy representation (i.e., the linear isotropy 
action) of a pseudo-Riemannian symmetric space (see Theorem 1 of [6]).  
Also, they showed that the principal orbits of the polar representation 
through a semi-simple element (i.e., the orbit through a regular element 
(in the sense of [6])) is an isoparametric submanifold in the sense of [6] 
by investigating the complexified representation (see Theorem 11 
(also Example 12) of [6]), 
where an isoparametric submanifold in the sense of [6] means 
the finite dimensional version of a proper complex isoparametric submanifold 
in a pseudo-Hilbert space defined in [14].  
All isoparametric submanifold (in a pseudo-Euclidean space) in this sense 
are isoparametric manifolds with flat section in the sense of [11].   
On the other hand, we [20] showed that, for a Hermann type action 
$H\curvearrowright G/K$ (i.e., $H$ is a symmetric subgroup of $G$) on a 
(semi-simple) pseudo-Riemannian symmetric space $G/K$, the principal 
$H$-orbits through $\exp_G(w)K$ ($w:$a semi-simple element, 
$\exp_G:$the exponential map of $G$) are curvature-adapted proper complex 
equifocal submanifolds (hence isoparametric submanifolds with flat section 
in the sense of [11]).  

In this paper, we shall first define the complexification $f^{\Bbb C}$ of 
a $C^{\omega}$-map of a $C^{\omega}$-Koszul manifold $(M,\nabla)$ into another 
$C^{\omega}$-Koszul manifold $(\widetilde M,\widetilde{\nabla})$ as a map of 
a tubular neighborhood $(M_{\nabla}^{\Bbb C})_f$ of $M$ 
in $M_{\nabla}^{\Bbb C}$ into ${\widetilde M}^{\Bbb C}_{\widetilde{\nabla}}$ 
and show that $f^{\Bbb C}$ is holomorphic and that, if $f$ is an immersion, 
then $f^{\Bbb C}$ also is an immersion on a tubular neighborhood 
$(M_{\nabla}^{\Bbb C})_{f,i}$ of $M$ in $(M_{\nabla}^{\Bbb C})_f$ 
(see Section 4).  
Let $(M,g)$ be a $C^{\omega}$-pseudo-Riemannian manifold.  Next, on a 
tubular neighborhood $(M_g^{\Bbb C})_A$ (which we take as largely as 
possible) of $M$ in $M_g^{\Bbb C}$, we define an anti-Kaehler metric $g_A$ 
compatible with $J^g$ 
(i.e., $g_A(J^gX,J^gY)=-g_A(X,Y)\,\,(X,Y\in TU^g_A),\,\,\,\nabla J^g=0$) 
satisfying $\iota^{\ast}g_A=g$, where $\nabla$ is the Levi-Civita connection 
of $g_A$ and $\iota$ is the inclusion map of $M$ into $(M_g^{\Bbb C})_A$.  
Note that $g_A$ is defined uniquely.  
We show that, for a $C^{\omega}$-isometric immersion 
$f:(M,g)\hookrightarrow(\widetilde M,\widetilde g)$ between 
$C^{\omega}$-pseudo-Riemannian manifolds, 
$f^{\Bbb C}:((M_g^{\Bbb C})_A\cap (M_g^{\Bbb C})_{f:i},\,\,g_A)\to
((\widetilde M_{\widetilde g}^{\Bbb C})_A,\,\,\widetilde g_A)$ 
is a holomorphic and isometric (that is, an anti-Kaehler) immersion.  
Next, for a pseudo-Riemannian homogeneous space, we define its another 
complexification as the quotient of the complexification of its isometry 
group by the complexification of its isotropy group, where we assume that the 
isometry group and the isotropy group have faithful real representations.  
Note that this quotient has a natural anti-Kaehler structure.  
The first purpose of this paper is to find an interesting relation between 
two complexifications (see Theorem 6.1).  
The second purpose of this paper is to define the dual of 
a $C^{\omega}$-pseudo-Riemannian manifold $(M,g)$ at each point and the dual 
of a totally geodesic $C^{\omega}$-submanifold of $(M,g)$ in 
the anti-Kaehler manifold $((M_g^{\Bbb C})_A,g_A)$ (see Definitions 2 and 
3 in Section 7).  
Next we define the notions of a complex Jacobi field in an anti-Kaehler 
manifold and a complex focal radius of an anti-Kaehler submanifold 
and show some facts related to them (see Section 8).  
Furthermore, we define the notions of 
a complex equifocal submanifold  and an isoparametric one in a 
pseudo-Riemannian homogeneous space and investigate the equivalence between 
their notions for a $C^{\omega}$-submanifold in a pseudo-Riemannian symmetric 
space (see Section 9).  
The third purpose of this paper is to show that, almost all orbits of 
the $G$-action on the complexification $(((G/K)_g^{\Bbb C})_A,g_A)$ of a 
pseudo-Riemannian symmetric space $(G/K,g)$ are curvature-adapted 
isoparametric submanifolds with flat section such that the shape operators 
are complex diagonalizable (see Theorem 9.3).  

\vspace{0.5truecm}

\noindent
{\bf Future plan of research.} {\sl 
We plan to solve both of various problems (for example, 
problems for harmonic analysis) in a $C^{\omega}$-pseudo-Riemannian 
manifold $(M,g)$ and the corresponding problems in the dual of $(M,g)$ 
by solving the corresponding problems in $((M_g^{\Bbb C})_A,g_A)$.}

\section{Basic notions and facts} 
In this section, we shall recall basic notions and facts.  Let $(M,\nabla)$ 
be a $C^{\infty}$-Koszul manifold and $\pi:TM\to M$ be the tangent bundle 
of $M$.  Denote by $\displaystyle{\mathop{TM}^{\circ}}$ the punctured tangent 
bundle $TM\setminus M$, where $M$ is identified with the zero section of 
$TM$.  Denote by $\mathfrak V$ the vertical distribution on $TM$ and by 
$\mathfrak H$ the horizontal distribution on $TM$ with respect to $\nabla$.  
Also, denote by $w_u^V(\in\mathfrak V_u)$ the vertical lift of 
$w\in T_{\pi(u)}M$ to $u$.  Let $\Phi_t$ be the geodesic flow of $\nabla$ and 
$X^S$ be the vector field on $TM$ associated with $\Phi_t$.  Define a 
distribution $\mathfrak L^{\nabla}$ on $\displaystyle{\mathop{TM}^{\circ}}$ 
by $\mathfrak L^{\nabla}_u:={\rm Span}\{u^V_u,X^S_u\}$ 
($u\in\displaystyle{\mathop{TM}^{\circ}}$).  
This distribution  $\mathfrak L^{\nabla}$ is involutive and hence defines 
a foliation on $\displaystyle{\mathop{TM}^{\circ}}$.  This foliation is 
called the {\it Koszul foliation} and we denote it by $\mathfrak F^{\nabla}$.  
In particular, if $\nabla$ is the Levi-Civita connection of a 
pseudo-Riemannian metric, then we call it a {\it Levi-Civita foliation}.  
These terminologies are used in [30].  Let $\gamma:I\to M$ be a maximal 
geodesic.  The image $\gamma_{\ast}(\displaystyle{\mathop{TI}^{\circ})}$ 
yields two leaves of $\mathfrak F^{\nabla}$ and all leaves of 
$\mathfrak F^{\nabla}$ are obtained in this way.  Let $\xi$ be a vector field 
along $\gamma_{\ast}$.  If there exists a geodesic variation $\gamma_t$ in 
$M$ satisfying $\gamma_0=\gamma$ and $\frac{d}{dt}\vert_{t=0}\gamma_{t\ast}
=\xi$, then $\xi$ is called a {\it parallel vector field}.  Note that $\xi$ is 
an extension of the Jacobi field $\frac{d}{dt}\vert_{t=0}\gamma_t$ along 
$\gamma$.  If $(M,\nabla)$ is a $C^{\omega}$-Koszul manifold, then there 
uniquely exists a complex structure $J^{\nabla}$ on a suitable domain 
$U^{\nabla}$ of $TM$ containing $M$ such that for each maximal geodesic 
$\gamma$ in $(M,\nabla)$, $\gamma_{\ast}:\gamma_{\ast}^{-1}(U^{\nabla})\to
(U^{\nabla},J^{\nabla})$ is holomorphic (see Theorem 0.3 of [30]), 
where $\gamma_{\ast}^{-1}(U^{\nabla})$ is regarded as an open set of 
${\Bbb C}$ under the natural identification of $T{\Bbb R}$ with ${\Bbb C}$.  
We take $U^{\nabla}$ as largely as possible.  This complex structure 
$J^{\nabla}$ is called the {\it adapted complex structure}.  We denote this 
complex manifold $(U^{\nabla},J^{\nabla})$ by $M_{\nabla}^{\Bbb C}$ and call it 
the complexification of $(M,\nabla)$.  In particular, if $\nabla$ is the 
Levi-Civita connection of a pseudo-Riemannian metric $g$, then 
$U^{\nabla},J^{\nabla}$ and $M_{\nabla}^{\Bbb C}$ are denoted by $U^g,\,J^g$ 
and $M_g^{\Bbb C}$, respectively.  Denote by $R$ the curvature tensor of 
$\nabla$.  According to Remark 2.2 of [5] and the statement (b) of Page 8 of 
[5], we see that, if $(M,\nabla)$ is locally symmetric (i.e., $\nabla\,:\,$ 
torsion-free and $\nabla R=0$) and the spectrum of $R(\cdot,X)X$ contains no 
negative number for each $X\in TM$, then the adapted complex structure 
$J^{\nabla}$ is defined on $TM$ (i.e., $U^{\nabla}=TM$).  

\section{Anti-Kaehler manifolds}
Let $M$ be a $C^{\infty}$-manifold, $J$ be a complex structure on $M$ and $g$ 
be a pseudo-Riemannian metric on $M$.  Denote by $\nabla$ the Levi-Civita 
connection of $g$.  If $g(JX,JY)=-g(X,Y)$ for any tangent vectors $X$ and $Y$ 
of $M$, then $(M,J,g)$ is called a {\it anti-Hermitian manifold}.  
Furthermore, if $\nabla J=0$, then it is called an 
{\it anti-Kaehler manifold}.  
For an anti-Kaehler manifold, the following remarkable fact holds.  

\vspace{0.5truecm}

\noindent
{\bf Proposition 3.1.} {\sl Let $(M,J,g)$ be an anti-Kaehler manifold and 
$\exp_p$ be the exponential map at $p\in M$.  Then $\exp_p:(T_pM,J_p)\to
(M,J)$ is holomorphic.}

\vspace{0.5truecm}

\noindent
{\it Proof.} Let $u\in T_pM$ and $X\in T_u(T_pM)$.  Define a geodesic 
variation $\delta$ (resp. $\bar{\delta}$) by $\delta(t,s):=\exp_p(t(u+sX))$ 
(resp. $\bar{\delta}(t,s):=\exp_p(t(u+sJ_pX)))$ for $(t,s)\in[0,1]^2$.  
Let $Y:=\delta_{\ast}(\frac{\partial}{\partial s}\vert_{s=0})$ and 
$\bar Y:=\bar{\delta}_{\ast}
(\frac{\partial}{\partial s}\vert_{s=0})$, which are Jacobi fields along 
the geodesic $\gamma_u$ with $\dot{\gamma}_u(0)=u$.  Since $(M,J,g)$ is 
anti-Kaehler, we have $\nabla J=0$ and $R(Jv,w)=JR(v,w)$ ($v,w\in TM$) 
(by Lemma 5.2 of [1]), where 
$R$ is the curvature tensor of $g$.  Hence we have 
$$\nabla_{\dot{\gamma}_u}\nabla_{\dot{\gamma}_u}(JY)+R(JY,\dot{\gamma}_u)
\dot{\gamma}_u=J\left(\nabla_{\dot{\gamma}_u}\nabla_{\dot{\gamma}_u}Y
+R(Y,\dot{\gamma}_u)\dot{\gamma}_u\right)=0,$$
that is, $JY$ is also a Jacobi field along $\gamma_u$.  Also, we have 
$JY(0)=\bar Y(0)=0$ and $\nabla_{\dot{\gamma}_u(0)}JY=
\nabla_{\dot{\gamma}_u(0)}\bar Y=J_pX$.  Hence we have $JY=\bar Y$.  
On the other hand, we have $JY(1)=J_{\gamma_u(1)}(\exp_p)_{\ast u}(X)$ and 
$\bar Y(1)=(\exp_p)_{\ast u}(J_pX)$.  Therefore $J_{\gamma_u(1)}\circ
(\exp_p)_{\ast u}=(\exp_p)_{\ast u}\circ J_p$ follows from the arbitrariness 
of $X$.  Since this relation holds for any $u\in T_pM$, 
$\exp_p\,:\,(T_pM,J_p)\to (M,J)$ is holomorphic. \hspace{6.4truecm}q.e.d.

\vspace{0.5truecm}

According to this fact, we can define so-called normal holomorphic coordinate 
around each point $p$ of a real $2n$-dimensional anti-Kaehler manifold 
$(M,J,g)$ as follows.  
Let $\widetilde U$ be a neighborhood of the origin of $T_pM$ such that 
$\exp_p\vert_{\widetilde U}$ is a diffeomorphism and $(e_1,J_pe_1,\cdots,e_n,
J_pe_n)$ be a $J_p$-base of $T_pM$.  Define $\widetilde{\phi}:{\Bbb C}^n\to 
T_pM$ by $\widetilde{\phi}(x_1+\sqrt{-1}y_1,\cdots,$\newline
$x_n+\sqrt{-1}y_n)=\sum\limits_{i=1}^n(x_ie_i+y_iJ_pe_i)$.  Set $U:=\exp_p(\widetilde U)$ and 
$\phi:=\widetilde{\phi}^{-1}\circ(\exp_p\vert_{\widetilde U})^{-1}$.  
According to Proposition 3.1, $(U,\phi)$ is a holomorphic local coordinate of 
$(M,J,g)$.  We call such a coordinate a {\it normal holomorphic coordinate} 
of $(M,J,g)$.  
Let $v\in T_pM$ and define a map $\gamma_v^{\Bbb C}:D\to M$ by 
$\gamma_v^{\Bbb C}(z)=\exp_p(({\rm Re}\,z)v+({\rm Im}\,z)J_pv)$ ($z\in D$), 
where $D$ is an open neighborhood of $0$ in ${\Bbb C}$.  We may assume that 
$\gamma_v^{\Bbb C}$ is an immersion by shrinking $D$ if necessary.  
According to Proposition 3.1, $\gamma_v^{\Bbb C}$ is the holomorphic 
extension of $\gamma_v$ and hence it is totally geodesic.  
We call $\gamma_v^{\Bbb C}$ a {\it complex geodesic} in $(M,J,g)$.  

Next we give examples of 
an anti-Kaehler manifold.  Let $(G,K)$ be a semi-simple symmetric pair 
and $\mathfrak g=\mathfrak k+\mathfrak p$ be the canonical 
decomposition of $\mathfrak g:={\rm Lie}\,G$ associated with $(G,K)$.  
Denote by $g$ the $G$-invariant pseudo-Riemannian metric on 
a quotient manifold $G/K$ arising from the restriction 
$B\vert_{\mathfrak p\times\mathfrak p}$ to 
$\mathfrak p$ of the Killing form $B$ of $\mathfrak g$.  
Then $(G/K,g)$ and $(G/K,-g)$ are (semi-simple) pseudo-Riemannian symmetric 
spaces.  Note that $(G/K,-g)$ is a Riemannian symmetric space of compact type 
if $(G,K)$ is a Riemannian symmetric pair of compact type and that $(G/K,g)$ 
is a Riemannian symmetric space of non-compact type if $(G,K)$ is a Riemannian 
symmetric pair of non-compact type.  
Let $G^{\Bbb C},K^{\Bbb C},
\mathfrak g^{\Bbb C},\mathfrak k^{\Bbb C}$ and $\mathfrak p^{\Bbb C}$ be the 
complexifications of $G, K, \mathfrak g, \mathfrak k$ and $\mathfrak p$, 
respectively.  For the complexification $B^{\Bbb C}\,(:\,\mathfrak g^{\Bbb C}
\times\mathfrak g^{\Bbb C}\to{\Bbb C})$ of $B$, $2{\rm Re}\,B^{\Bbb C}$ is 
the Killing form of $\mathfrak g^{\Bbb C}$ regarded as a real Lie algebra, 
where ${\rm Re}\,B^{\Bbb C}$ is the real part of $B^{\Bbb C}$.  The pair 
$(G^{\Bbb C},K^{\Bbb C})$ is a semi-simple symmetric pair, where $G^{\Bbb C}$ and 
$K^{\Bbb C}$ are regarded as real Lie groups.  
Denote by $\widetilde g$ the $G^{\Bbb C}$-invariant pseudo-Riemannian metric 
on $G^{\Bbb C}/K^{\Bbb C}$ arising from $2{\rm Re}\,B^{\Bbb C}
\vert_{\mathfrak p^{\Bbb C}\times\mathfrak p^{\Bbb C}}$ and 
by $J$ the $G^{\Bbb C}$-invariant complex structure arising from 
$j:\mathfrak p^{\Bbb C}\to\mathfrak p^{\Bbb C}\,\,
(\displaystyle{\mathop{\Leftrightarrow}_{{\rm def}}\,jX=\sqrt{-1}X)}$.  Then 
$(G^{\Bbb C}/K^{\Bbb C},J,\widetilde g)$ and 
$(G^{\Bbb C}/K^{\Bbb C},J,-\widetilde g)$ are anti-Kaehler manifolds.  We call 
these anti-Kaehler manifolds 
the {\it anti-Kaehler symmetric spaces associated with} $(G/K,g)$ 
and $(G/K,-g)$, respectively.  See [19] about general theory of 
an anti-Kaehler symmetric space.  

\section{A complexification of a $C^{\omega}$-map between Koszul manifolds}
In this section, we shall define the complexification of a $C^{\omega}$-map 
between $C^{\omega}$-Koszul manifolds and investigate it.  
Let $f\,:\,(M,\nabla)\to(\widetilde M,\widetilde{\nabla})$ be 
a $C^{\omega}$-map between $C^{\omega}$-Koszul manifolds.  First we shall 
recall the definition of the (maximal) holomorphic extension $\alpha^h$ of a 
$C^{\omega}$-curve 
$\alpha:(a,b)\to\widetilde M$ in $\widetilde M_{\widetilde{\nabla}}^{\Bbb C}$.  
Fix $t_0\in(a,b)$ and take a holomorphic local coordinate 
$(V,\phi=(z_1,\cdots,z_m))$ of $\widetilde M_{\widetilde{\nabla}}^{\Bbb C}$ 
around $\alpha(t_0)$ satisfying $\widetilde M\cap V=\phi^{-1}({\Bbb R}^m)$, 
where $m={\rm dim}\,\widetilde M$.  Let $(\phi\circ\alpha)(t)=
(\alpha_1(t),\cdots,\alpha_m(t))$.  Since $\alpha_i(t)$ ($i=1,\cdots,m$) are of class $C^{\omega}$, we get their holomorphic extensions $\alpha_i^h:D_i\to
{\Bbb C}$ ($i=1,\cdots,m$), where $D_i$ is a neighborhood of $t_0$ in 
${\Bbb C}$.  Define $\alpha_{t_0}^h\,:\,\displaystyle{
\left(\mathop{\cap}_{i=1}^k D_i\right)\cap(\alpha_1^h\times\cdots\times
\alpha_m^h)^{-1}(\phi(V))\to M_{\nabla}^{\Bbb C}}$ by $\alpha_{t_0}^h(z):=
\phi^{-1}(\alpha_1^h(z),\cdots,\alpha_m^h(z))$.  This complex curve 
$\alpha_{t_0}^h$ is a holomorphic extension of 
$\alpha\vert_{(t_0-\varepsilon,t_0+\varepsilon)}$, where $\varepsilon$ is 
a sufficiently small positive number.  For each $t\in(a,b)$, we get a 
holomorphic extension $\alpha_t^h$ of 
$\alpha\vert_{(t-\varepsilon',t+\varepsilon')}$, where $\varepsilon'$ is 
a sufficiently small positive number.  
By patching $\{\alpha_t^h\}_{t\in(a,b)}$, we get a holomorphic extension of 
$\alpha$ and furthermore, by extending the holomorphic extension to the 
maximal one, we get the maximal holomorphic extension $\alpha^h$.  
Now we shall define the complexification 
$f^{\Bbb C}$ of $f$.  

\vspace{0.5truecm}

\noindent
{\bf Definition.} Let $(M_{\nabla}^{\Bbb C})_f:=\{v\in 
M_{\nabla}^{\Bbb C}\,\vert\,\sqrt{-1}\in{\rm Dom}((f\circ\gamma_v)^h)\}$, 
where $\gamma_v$ is 
the geodesic in $(M,\nabla)$ with $\dot{\gamma}_v(0)=v$, $(f\circ\gamma_v)^h$ 
is the (maximal) holomorphic extension of $f\circ\gamma_v$ in 
$\widetilde M_{\widetilde{\nabla}}^{\Bbb C}$ and ${\rm Dom}((f\circ\gamma_v)^h)$ is the domain of $(f\circ\gamma_v)^h$.  This set $(M_{\nabla}^{\Bbb C})_f$ 
is a tubular neighborhood of $M$ in $M_{\nabla}^{\Bbb C}$.  
We define $f^{\Bbb C}:(M_{\nabla}^{\Bbb C})_f\to
\widetilde M_{\widetilde{\nabla}}^{\Bbb C}$ by $f^{\Bbb C}(v):=
(f\circ\gamma_v)^h(\sqrt{-1})$ ($v\in(M_{\nabla}^{\Bbb C})_f$).  

\vspace{0.3truecm}

For this complexification $f^{\Bbb C}$, the following facts hold.  

\vspace{0.3truecm}

\noindent
{\bf Proposition 4.1.} {\sl Let $f:(M,\nabla)\to
(\widetilde M,\widetilde{\nabla})$ 
be a $C^{\omega}$-map between $C^{\omega}$-Koszul manifolds.  
Then $f^{\Bbb C}:(M_{\nabla}^{\Bbb C})_f\to
\widetilde M_{\widetilde{\nabla}}^{\Bbb C}$ 
is the (maximal) holomorphic extension of $f$.  
Also, if $f$ is an immersion, then $f^{\Bbb C}$ is an 
immersion on a tubular neighborhood 
(which is denoted by $(M_{\nabla}^{\Bbb C})_{f:i}$ in the sequel) of $M$ in 
$(M_{\nabla}^{\Bbb C})_f$.}

\vspace{0.3truecm}

\centerline{
\unitlength 0.1in
\begin{picture}( 61.0700, 27.0900)(  8.3200,-33.8900)
%
\special{pn 8}%
\special{pa 3204 1628}%
\special{pa 3204 2684}%
\special{fp}%
\special{pa 3204 2684}%
\special{pa 4134 2684}%
\special{fp}%
\special{pa 4134 2684}%
\special{pa 4134 1628}%
\special{fp}%
\special{pa 4134 1628}%
\special{pa 3204 1628}%
\special{fp}%
%
\special{pn 8}%
\special{pa 3204 2156}%
\special{pa 4134 2156}%
\special{fp}%
%
\special{pn 8}%
\special{pa 3204 1760}%
\special{pa 4134 1760}%
\special{pa 4134 2552}%
\special{pa 3204 2552}%
\special{pa 3204 1760}%
\special{dt 0.045}%
%
\special{pn 8}%
\special{pa 3204 1760}%
\special{pa 4134 1760}%
\special{pa 4134 2156}%
\special{pa 3204 2156}%
\special{pa 3204 1760}%
\special{dt 0.045}%
%
\special{pn 20}%
\special{pa 3442 2148}%
\special{pa 3680 2148}%
\special{dt 0.054}%
\special{sh 1}%
\special{pa 3680 2148}%
\special{pa 3612 2128}%
\special{pa 3626 2148}%
\special{pa 3612 2168}%
\special{pa 3680 2148}%
\special{fp}%
%
\special{pn 20}%
\special{sh 1}%
\special{ar 3442 1924 10 10 0  6.28318530717959E+0000}%
\special{sh 1}%
\special{ar 3442 1924 10 10 0  6.28318530717959E+0000}%
%
\special{pn 8}%
\special{pa 3442 1628}%
\special{pa 3442 2684}%
\special{fp}%
%
\special{pn 13}%
\special{pa 3442 2156}%
\special{pa 3982 2148}%
\special{fp}%
\put(34.0800,-18.8400){\makebox(0,0)[rt]{$v$}}%
\put(37.1800,-22.0100){\makebox(0,0)[rt]{$v$}}%
%
\special{pn 8}%
\special{pa 4042 2756}%
\special{pa 3864 2162}%
\special{dt 0.045}%
\special{sh 1}%
\special{pa 3864 2162}%
\special{pa 3864 2232}%
\special{pa 3880 2212}%
\special{pa 3902 2220}%
\special{pa 3864 2162}%
\special{fp}%
%
\special{pn 8}%
\special{pa 3870 1514}%
\special{pa 4076 2156}%
\special{dt 0.045}%
\special{sh 1}%
\special{pa 4076 2156}%
\special{pa 4074 2086}%
\special{pa 4060 2104}%
\special{pa 4036 2098}%
\special{pa 4076 2156}%
\special{fp}%
%
\special{pn 8}%
\special{pa 3132 1396}%
\special{pa 3442 1706}%
\special{dt 0.045}%
\special{sh 1}%
\special{pa 3442 1706}%
\special{pa 3408 1646}%
\special{pa 3404 1668}%
\special{pa 3380 1674}%
\special{pa 3442 1706}%
\special{fp}%
\put(41.6700,-13.6900){\makebox(0,0)[rt]{$0$-section$=M$}}%
\put(32.4300,-12.5000){\makebox(0,0)[rt]{$T_{\gamma_v(0)}M$}}%
\put(41.0800,-27.7500){\makebox(0,0)[rt]{$\gamma_v$}}%
\put(38.1700,-31.1200){\makebox(0,0)[rt]{$TM$}}%
%
\special{pn 8}%
\special{ar 3172 2156 152 396  1.5707963 1.6145919}%
\special{ar 3172 2156 152 396  1.7459788 1.7897744}%
\special{ar 3172 2156 152 396  1.9211613 1.9649569}%
\special{ar 3172 2156 152 396  2.0963438 2.1401394}%
\special{ar 3172 2156 152 396  2.2715263 2.3153219}%
\special{ar 3172 2156 152 396  2.4467087 2.4905044}%
\special{ar 3172 2156 152 396  2.6218912 2.6656868}%
\special{ar 3172 2156 152 396  2.7970737 2.8408693}%
\special{ar 3172 2156 152 396  2.9722562 3.0160518}%
\special{ar 3172 2156 152 396  3.1474387 3.1912343}%
\special{ar 3172 2156 152 396  3.3226211 3.3664168}%
\special{ar 3172 2156 152 396  3.4978036 3.5415992}%
\special{ar 3172 2156 152 396  3.6729861 3.7167817}%
\special{ar 3172 2156 152 396  3.8481686 3.8919642}%
\special{ar 3172 2156 152 396  4.0233511 4.0671467}%
\special{ar 3172 2156 152 396  4.1985336 4.2423292}%
\special{ar 3172 2156 152 396  4.3737160 4.4175117}%
\special{ar 3172 2156 152 396  4.5488985 4.5926941}%
\put(29.9200,-20.2300){\makebox(0,0)[rt]{$(M^{\Bbb C}_{\nabla})_f$}}%
\put(46.6000,-21.6000){\makebox(0,0)[rt]{$\longrightarrow$ }}%
%
\special{pn 8}%
\special{ar 5884 1192 1056 264  3.1640428 6.2831853}%
\special{ar 5884 1192 1056 264  0.0000000 3.1415927}%
%
\special{pn 8}%
\special{ar 5884 3040 1056 264  3.1611432 6.2831853}%
\special{ar 5884 3040 1056 264  0.0000000 3.1415927}%
%
\special{pn 8}%
\special{pa 4828 1192}%
\special{pa 4828 3046}%
\special{fp}%
%
\special{pn 8}%
\special{pa 6940 1192}%
\special{pa 6940 3040}%
\special{fp}%
%
\special{pn 8}%
\special{ar 5884 2116 1056 264  6.2831853 6.2831853}%
\special{ar 5884 2116 1056 264  0.0000000 3.1415927}%
%
\special{pn 8}%
\special{ar 6408 2250 1188 264  3.1415927 4.7103475}%
%
\special{pn 8}%
\special{ar 6420 2770 1188 264  3.1415927 3.2242373}%
\special{ar 6420 2770 1188 264  3.2738241 3.3564687}%
\special{ar 6420 2770 1188 264  3.4060555 3.4887001}%
\special{ar 6420 2770 1188 264  3.5382869 3.6209315}%
\special{ar 6420 2770 1188 264  3.6705183 3.7531629}%
\special{ar 6420 2770 1188 264  3.8027497 3.8853943}%
\special{ar 6420 2770 1188 264  3.9349811 4.0176257}%
\special{ar 6420 2770 1188 264  4.0672125 4.1498571}%
\special{ar 6420 2770 1188 264  4.1994439 4.2820885}%
\special{ar 6420 2770 1188 264  4.3316753 4.4143199}%
\special{ar 6420 2770 1188 264  4.4639067 4.5465513}%
\special{ar 6420 2770 1188 264  4.5961381 4.6787827}%
%
\special{pn 8}%
\special{ar 6410 1760 1188 264  3.1415927 3.2242373}%
\special{ar 6410 1760 1188 264  3.2738241 3.3564687}%
\special{ar 6410 1760 1188 264  3.4060555 3.4887001}%
\special{ar 6410 1760 1188 264  3.5382869 3.6209315}%
\special{ar 6410 1760 1188 264  3.6705183 3.7531629}%
\special{ar 6410 1760 1188 264  3.8027497 3.8853943}%
\special{ar 6410 1760 1188 264  3.9349811 4.0176257}%
\special{ar 6410 1760 1188 264  4.0672125 4.1498571}%
\special{ar 6410 1760 1188 264  4.1994439 4.2820885}%
\special{ar 6410 1760 1188 264  4.3316753 4.4143199}%
\special{ar 6410 1760 1188 264  4.4639067 4.5465513}%
\special{ar 6410 1760 1188 264  4.5961381 4.6787827}%
%
\special{pn 8}%
\special{pa 5220 1740}%
\special{pa 5220 2308}%
\special{fp}%
%
\special{pn 8}%
\special{pa 5224 2366}%
\special{pa 5224 2776}%
\special{fp}%
%
\special{pn 8}%
\special{pa 6420 1520}%
\special{pa 6420 2372}%
\special{fp}%
%
\special{pn 8}%
\special{pa 6418 2380}%
\special{pa 6418 2512}%
\special{fp}%
%
\special{pn 8}%
\special{ar 5884 2116 1050 264  5.2998336 6.2831853}%
%
\special{pn 8}%
\special{ar 5882 2116 1056 272  4.1004459 5.1428855}%
%
\special{pn 8}%
\special{ar 5884 2116 1056 264  3.1415927 3.9399397}%
%
\special{pn 8}%
\special{pa 5488 1720}%
\special{pa 5488 2248}%
\special{fp}%
%
\special{pn 8}%
\special{ar 6274 2326 932 146  3.6646355 4.5525377}%
%
\special{pn 8}%
\special{pa 6134 2188}%
\special{pa 6134 1640}%
\special{fp}%
%
\special{pn 8}%
\special{ar 6312 1812 1050 172  3.7882236 4.5464689}%
%
\special{pn 4}%
\special{pa 3204 1760}%
\special{pa 4134 1760}%
\special{pa 4134 2552}%
\special{pa 3204 2552}%
\special{pa 3204 1760}%
\special{ip}%
%
\special{pn 4}%
\special{pa 3200 1760}%
\special{pa 4130 1760}%
\special{pa 4130 2552}%
\special{pa 3200 2552}%
\special{pa 3200 1760}%
\special{ip}%
%
\special{pn 4}%
\special{pa 5488 1720}%
\special{pa 6134 1720}%
\special{pa 6134 1654}%
\special{pa 5488 1654}%
\special{pa 5488 1720}%
\special{ip}%
%
\special{pn 20}%
\special{sh 1}%
\special{ar 5488 1832 10 10 0  6.28318530717959E+0000}%
\special{sh 1}%
\special{ar 5488 1832 10 10 0  6.28318530717959E+0000}%
%
\special{pn 13}%
\special{ar 6416 2250 1188 258  3.7912491 4.4632428}%
\put(45.5000,-20.1000){\makebox(0,0)[rt]{$f^{\Bbb C}$}}%
%
\special{pn 4}%
\special{pa 5210 1742}%
\special{pa 6406 1490}%
\special{pa 6406 2526}%
\special{pa 6062 2534}%
\special{pa 6062 2534}%
\special{pa 5812 2546}%
\special{pa 5812 2546}%
\special{pa 5514 2606}%
\special{pa 5304 2678}%
\special{pa 5244 2732}%
\special{pa 5230 2758}%
\special{pa 5218 2778}%
\special{pa 5218 1748}%
\special{pa 5218 1748}%
\special{pa 5210 1742}%
\special{ip}%
%
\special{pn 8}%
\special{pa 5090 930}%
\special{pa 5816 1624}%
\special{dt 0.045}%
\special{sh 1}%
\special{pa 5816 1624}%
\special{pa 5782 1564}%
\special{pa 5778 1586}%
\special{pa 5754 1592}%
\special{pa 5816 1624}%
\special{fp}%
\special{pa 5500 1624}%
\special{pa 5500 1624}%
\special{dt 0.045}%
%
\special{pn 8}%
\special{pa 4756 1238}%
\special{pa 5600 1778}%
\special{dt 0.045}%
\special{sh 1}%
\special{pa 5600 1778}%
\special{pa 5554 1726}%
\special{pa 5554 1750}%
\special{pa 5532 1760}%
\special{pa 5600 1778}%
\special{fp}%
%
\special{pn 8}%
\special{pa 6372 2590}%
\special{pa 5962 2024}%
\special{dt 0.045}%
\special{sh 1}%
\special{pa 5962 2024}%
\special{pa 5986 2090}%
\special{pa 5994 2066}%
\special{pa 6018 2066}%
\special{pa 5962 2024}%
\special{fp}%
%
\special{pn 8}%
\special{pa 6068 1204}%
\special{pa 6320 1990}%
\special{dt 0.045}%
\special{sh 1}%
\special{pa 6320 1990}%
\special{pa 6318 1920}%
\special{pa 6304 1940}%
\special{pa 6280 1934}%
\special{pa 6320 1990}%
\special{fp}%
\put(60.8800,-33.8900){\makebox(0,0)[rt]{$T\widetilde M$}}%
\put(61.0800,-10.7900){\makebox(0,0)[rt]{$f(M)$}}%
\put(66.0000,-26.3000){\makebox(0,0)[rt]{$f\circ\gamma_v$}}%
\put(52.2000,-7.2000){\makebox(0,0)[rt]{$f^{\Bbb C}(M^{\Bbb C})$}}%
\put(47.2800,-11.5800){\makebox(0,0)[rt]{$(f\circ\gamma_v)^h$}}%
\put(47.6000,-15.9000){\makebox(0,0)[rt]{$f^{\Bbb C}(v)$}}%
%
\special{pn 8}%
\special{pa 6398 882}%
\special{pa 6564 2076}%
\special{dt 0.045}%
\special{sh 1}%
\special{pa 6564 2076}%
\special{pa 6574 2006}%
\special{pa 6556 2022}%
\special{pa 6534 2012}%
\special{pa 6564 2076}%
\special{fp}%
\put(70.6000,-6.8000){\makebox(0,0)[rt]{$0$-section$=\widetilde M$}}%
%
\special{pn 8}%
\special{pa 4770 1660}%
\special{pa 5480 1830}%
\special{dt 0.045}%
\special{sh 1}%
\special{pa 5480 1830}%
\special{pa 5420 1796}%
\special{pa 5428 1818}%
\special{pa 5412 1834}%
\special{pa 5480 1830}%
\special{fp}%
%
\special{pn 4}%
\special{pa 3200 1770}%
\special{pa 3210 1760}%
\special{fp}%
\special{pa 3200 1830}%
\special{pa 3270 1760}%
\special{fp}%
\special{pa 3200 1890}%
\special{pa 3330 1760}%
\special{fp}%
\special{pa 3200 1950}%
\special{pa 3390 1760}%
\special{fp}%
\special{pa 3200 2010}%
\special{pa 3450 1760}%
\special{fp}%
\special{pa 3200 2070}%
\special{pa 3510 1760}%
\special{fp}%
\special{pa 3200 2130}%
\special{pa 3570 1760}%
\special{fp}%
\special{pa 3200 2190}%
\special{pa 3630 1760}%
\special{fp}%
\special{pa 3200 2250}%
\special{pa 3690 1760}%
\special{fp}%
\special{pa 3200 2310}%
\special{pa 3750 1760}%
\special{fp}%
\special{pa 3200 2370}%
\special{pa 3810 1760}%
\special{fp}%
\special{pa 3200 2430}%
\special{pa 3870 1760}%
\special{fp}%
\special{pa 3200 2490}%
\special{pa 3930 1760}%
\special{fp}%
\special{pa 3200 2550}%
\special{pa 3990 1760}%
\special{fp}%
\special{pa 3260 2550}%
\special{pa 4050 1760}%
\special{fp}%
\special{pa 3320 2550}%
\special{pa 4110 1760}%
\special{fp}%
\special{pa 3380 2550}%
\special{pa 4130 1800}%
\special{fp}%
\special{pa 3440 2550}%
\special{pa 4130 1860}%
\special{fp}%
\special{pa 3500 2550}%
\special{pa 4130 1920}%
\special{fp}%
\special{pa 3560 2550}%
\special{pa 4130 1980}%
\special{fp}%
\special{pa 3620 2550}%
\special{pa 4130 2040}%
\special{fp}%
\special{pa 3680 2550}%
\special{pa 4130 2100}%
\special{fp}%
\special{pa 3740 2550}%
\special{pa 4130 2160}%
\special{fp}%
\special{pa 3800 2550}%
\special{pa 4130 2220}%
\special{fp}%
\special{pa 3860 2550}%
\special{pa 4130 2280}%
\special{fp}%
\special{pa 3920 2550}%
\special{pa 4130 2340}%
\special{fp}%
\special{pa 3980 2550}%
\special{pa 4130 2400}%
\special{fp}%
\special{pa 4040 2550}%
\special{pa 4130 2460}%
\special{fp}%
\special{pa 4100 2550}%
\special{pa 4130 2520}%
\special{fp}%
%
\special{pn 4}%
\special{pa 5480 1720}%
\special{pa 5490 2240}%
\special{pa 5590 2220}%
\special{pa 5800 2200}%
\special{pa 6130 2190}%
\special{pa 6130 1660}%
\special{pa 5880 1670}%
\special{pa 5730 1690}%
\special{pa 5500 1710}%
\special{pa 5500 1710}%
\special{pa 5480 1720}%
\special{ip}%
%
\special{pn 4}%
\special{pa 5498 2238}%
\special{pa 5490 2230}%
\special{fp}%
\special{pa 5548 2228}%
\special{pa 5490 2170}%
\special{fp}%
\special{pa 5600 2220}%
\special{pa 5488 2108}%
\special{fp}%
\special{pa 5654 2214}%
\special{pa 5486 2046}%
\special{fp}%
\special{pa 5710 2210}%
\special{pa 5486 1986}%
\special{fp}%
\special{pa 5764 2204}%
\special{pa 5484 1924}%
\special{fp}%
\special{pa 5820 2200}%
\special{pa 5484 1864}%
\special{fp}%
\special{pa 5878 2198}%
\special{pa 5482 1802}%
\special{fp}%
\special{pa 5936 2196}%
\special{pa 5480 1740}%
\special{fp}%
\special{pa 5994 2194}%
\special{pa 5510 1710}%
\special{fp}%
\special{pa 6052 2192}%
\special{pa 5564 1704}%
\special{fp}%
\special{pa 6112 2192}%
\special{pa 5620 1700}%
\special{fp}%
\special{pa 6130 2150}%
\special{pa 5676 1696}%
\special{fp}%
\special{pa 6130 2090}%
\special{pa 5730 1690}%
\special{fp}%
\special{pa 6130 2030}%
\special{pa 5784 1684}%
\special{fp}%
\special{pa 6130 1970}%
\special{pa 5836 1676}%
\special{fp}%
\special{pa 6130 1910}%
\special{pa 5890 1670}%
\special{fp}%
\special{pa 6130 1850}%
\special{pa 5948 1668}%
\special{fp}%
\special{pa 6130 1790}%
\special{pa 6006 1666}%
\special{fp}%
\special{pa 6130 1730}%
\special{pa 6064 1664}%
\special{fp}%
\special{pa 6130 1670}%
\special{pa 6120 1660}%
\special{fp}%
%
\special{pn 4}%
\special{pa 5220 1786}%
\special{pa 5260 1706}%
\special{fp}%
\special{pa 5220 1876}%
\special{pa 5332 1652}%
\special{fp}%
\special{pa 5220 1966}%
\special{pa 5398 1610}%
\special{fp}%
\special{pa 5220 2056}%
\special{pa 5452 1592}%
\special{fp}%
\special{pa 5220 2146}%
\special{pa 5506 1574}%
\special{fp}%
\special{pa 5220 2236}%
\special{pa 5558 1562}%
\special{fp}%
\special{pa 5220 2326}%
\special{pa 5606 1556}%
\special{fp}%
\special{pa 5220 2416}%
\special{pa 5654 1550}%
\special{fp}%
\special{pa 5220 2506}%
\special{pa 5702 1544}%
\special{fp}%
\special{pa 5220 2596}%
\special{pa 5752 1534}%
\special{fp}%
\special{pa 5220 2686}%
\special{pa 5802 1524}%
\special{fp}%
\special{pa 5252 2714}%
\special{pa 5850 1514}%
\special{fp}%
\special{pa 5318 2670}%
\special{pa 5900 1506}%
\special{fp}%
\special{pa 5380 2636}%
\special{pa 5946 1506}%
\special{fp}%
\special{pa 5436 2616}%
\special{pa 5990 1506}%
\special{fp}%
\special{pa 5488 2602}%
\special{pa 6036 1506}%
\special{fp}%
\special{pa 5538 2588}%
\special{pa 6080 1506}%
\special{fp}%
\special{pa 5588 2580}%
\special{pa 6126 1506}%
\special{fp}%
\special{pa 5636 2574}%
\special{pa 6170 1506}%
\special{fp}%
\special{pa 5686 2566}%
\special{pa 6216 1506}%
\special{fp}%
\special{pa 5734 2560}%
\special{pa 6260 1506}%
\special{fp}%
\special{pa 5782 2554}%
\special{pa 6306 1506}%
\special{fp}%
\special{pa 5828 2548}%
\special{pa 6350 1506}%
\special{fp}%
\special{pa 5876 2544}%
\special{pa 6396 1506}%
\special{fp}%
\special{pa 5924 2540}%
\special{pa 6410 1566}%
\special{fp}%
\special{pa 5970 2536}%
\special{pa 6410 1656}%
\special{fp}%
\special{pa 6016 2534}%
\special{pa 6410 1746}%
\special{fp}%
\special{pa 6062 2530}%
\special{pa 6410 1836}%
\special{fp}%
\special{pa 6110 2528}%
\special{pa 6410 1926}%
\special{fp}%
\special{pa 6156 2526}%
\special{pa 6410 2016}%
\special{fp}%
%
\special{pn 4}%
\special{pa 6202 2524}%
\special{pa 6410 2106}%
\special{fp}%
\special{pa 6248 2520}%
\special{pa 6410 2196}%
\special{fp}%
\special{pa 6294 2518}%
\special{pa 6410 2286}%
\special{fp}%
\special{pa 6340 2516}%
\special{pa 6410 2376}%
\special{fp}%
\special{pa 6386 2516}%
\special{pa 6410 2466}%
\special{fp}%
\end{picture}%
\hspace{5.5truecm}}

\vspace{0.3truecm}

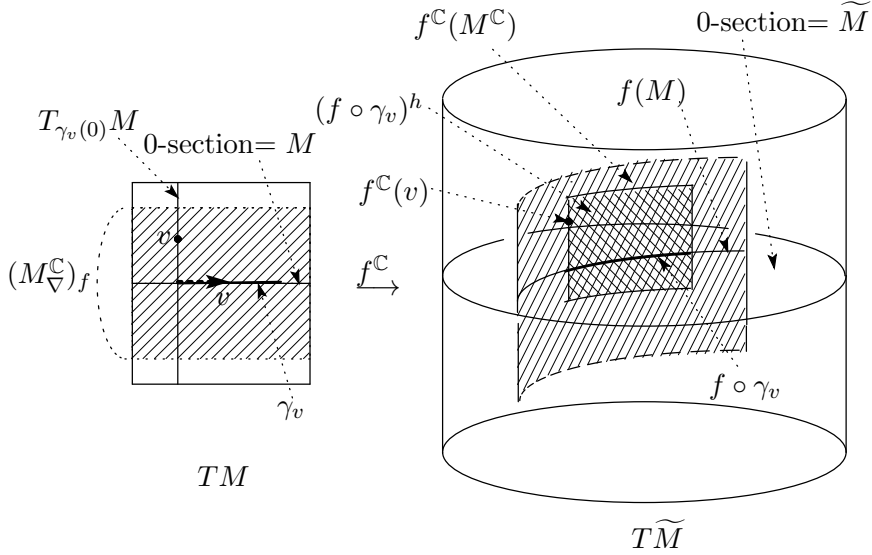
\centerline{{\bf Figure 1.}}

\vspace{0.5truecm}

\noindent
{\it Proof.} 
First we shall show $f^{\Bbb C}\vert_M=f$.  Take an arbitrary $p(=0_p)\in M\,
(=$the zero section of $TM$), where $0_p$ is the zero vector of $T_pM$.  
We have $f^{\Bbb C}(p)=f^{\Bbb C}(0_p)=(f\circ\gamma_{0_p})^h(\sqrt{-1})=f(p)$.  
Thus $f^{\Bbb C}\vert_M=f$ holds.  Next we shall show that $f^{\Bbb C}$ is 
holomorphic.  According to Theorem 3.4 of [27], we suffice to show that, 
for each geodesic $\gamma$ in $(M,\nabla)$, $f^{\Bbb C}\circ\gamma_{\ast}$ is 
holomorphic.  For each $z=x+\sqrt{-1}y\in{\rm Dom}
(f^{\Bbb C}\circ\gamma_{\ast})$, we have 
$$\begin{array}{l}
(f^{\Bbb C}\circ\gamma_{\ast})(z)=(f^{\Bbb C}\circ\gamma_{\ast})
(y(\frac{d}{dt})_x)=f^{\Bbb C}(y\dot{\gamma}(x))\\
\hspace{1.6truecm}=(f\circ\gamma_{y\dot{\gamma}(x)})^h(\sqrt{-1})
=(f\circ\gamma)^h(z),
\end{array}$$
where we note that the tangent bundle $T{\Bbb R}$ is identified with ${\Bbb C}$ 
under the correspondence $y(\frac{d}{dt})_x\leftrightarrow x+\sqrt{-1}y$.  
That is, we get $f^{\Bbb C}\circ\gamma_{\ast}=(f\circ\gamma)^h$.  Hence 
$f^{\Bbb C}\circ\gamma_{\ast}$ is holomorphic.  Thus the first-half part of 
the statement is shown.  The second-half part of the statement is trivial.  
\begin{flushright}q.e.d.\end{flushright}

\vspace{0.5truecm}

Let $(\widetilde M,\widetilde{\nabla})$ be an $m$-dimensional 
$C^{\omega}$-Koszul manifold, $F$ be a ${\Bbb R}^k$-valued $C^{\omega}$-
function over an open set $V$ of $\widetilde M$ ($k<m$) and $\bf a$ be a 
regular value of $F$.  Let $M:=F^{-1}({\bf a})$ and $\iota$ be the inclusion 
map of $M$ into $\widetilde M$.  Take an arbitrary $C^{\omega}$-Koszul 
connection $\nabla$ of $M$.  Then we have the following fact.  

\vspace{0.5truecm}

\noindent
{\bf Proposition 4.2.} {\sl The image $\iota^{\Bbb C}
((M_{\nabla}^{\Bbb C})_{\iota})$ is an open potion of $(F^h)^{-1}({\bf a})$, 
where $F^h$ is the (maximal) holomorphic extension of $F$ to 
$\widetilde M_{\widetilde{\nabla}}^{\Bbb C}$ (which is a ${\Bbb C}^k$-valued 
holomorphic function on a tubular neighborhood $\widetilde V$ of $V$ in 
$\widetilde M_{\widetilde{\nabla}}^{\Bbb C}$).}

\vspace{0.5truecm}

Here we shall explain the (maximal) holomorphic extension $F^h$ of $F$ to 
$\widetilde M^{\Bbb C}_{\widetilde{\nabla}}$.  Fix $p_0\in V$ and take 
a holomorphic local coordinate 
$(W_{p_0},\phi=(z_1,\cdots,z_m))$ of $\widetilde M^{\Bbb C}_{\widetilde{\nabla}}$ about $p_0$ satisfying $\widetilde M\cap W_{p_0}=\phi^{-1}({\Bbb R}^m)$ and 
$\widetilde M\cap W_{p_0}\subset V$.  
Since $F\circ(\phi\vert_{\widetilde M\cap W_{p_0}})^{-1}$ is of class 
$C^{\omega}$, we get its holomorphic extension 
$(F\circ(\phi\vert_{\widetilde M\cap W_{p_0}})^{-1})^h:D\to{\Bbb C}^k$, 
where $D$ is a neighborhood of $\phi(p_0)$ in 
${\Bbb C}^m$.  Define $F^h_{p_0}:\phi^{-1}(D\cap \phi(W_{p_0}))\to{\Bbb C}^k$ by 
$F^h_{p_0}:=(F\circ(\phi\vert_{\widetilde M\cap W_{p_0}})^{-1})^h\circ\phi
\vert_{\phi^{-1}(D\cap\phi(W_{p_0}))}$.  
This ${\Bbb C}^k$-valued function $F^h_{p_0}$ is 
a holomorphic extension of $F\vert_{\widetilde M\cap W_{p_0}}$ to 
$\widetilde M_{\widetilde{\nabla}}^{\Bbb C}$.  For each $p\in V$, we get a 
holomorphic extension $F^h_p$ of $F\vert_{V_p}$ ($V_p\,:\,$ a sufficiently 
small neighborhood of $p$ in $V$).  By patching $\{F^h_p\}_{p\in V}$, 
we get a holomorphic extension 
of $F$ and furthermore, by extending to the holomorphic extension to the 
maximal one, we get the desired (maximal) holomorphic extension $F^h$.  

\vspace{0.5truecm}

\noindent
{\it Proof of Proposition 4.2.} Take $X\in M_{\nabla}^{\Bbb C}\,(\subset TM)$ 
and $\gamma_X:(-\varepsilon,\varepsilon)\to M$ be the geodesic in $(M,\nabla)$ 
with $\dot{\gamma}_X(0)=X$.  Since $\gamma_X(t)\in M$, we have 
$F(\gamma_X(t))={\bf a}$, where $t\in(-\varepsilon,\varepsilon)$.  Let 
$(\iota\circ\gamma_X)^h(:D\to\widetilde M_{\widetilde{\nabla}}^{\Bbb C})$ be 
the (maximal) holomorphic extension of $\iota\circ\gamma_X$ in 
$\widetilde M_{\widetilde{\nabla}}^{\Bbb C}$.  Since 
$F^h\circ(\iota\circ\gamma_X)^h:((\iota\circ\gamma_X)^h)^{-1}(\widetilde V)\to
{\Bbb C}^k$ is holomorphic and $(F^h\circ(\iota\circ\gamma_X)^h)(t)={\bf a}$ 
($t\in(-\varepsilon,\varepsilon)$), we get $F^h\circ(\iota\circ\gamma_X)^h
\equiv{\bf a}$.  
Hence we get $F^h(\iota^{\Bbb C}(X))=F^h((\iota\circ\gamma_X)^h(\sqrt{-1}))=
{\bf a}$, that is, $\iota^{\Bbb C}(X)\in (F^h)^{-1}
({\bf a})$.  From the arbitrariness of $X$, it follows that 
$\iota^{\Bbb C}((M_{\nabla}^{\Bbb C})_{\iota})\subset (F^h)^{-1}({\bf a})$.  
Furthermore, since ${\rm dim}\,\iota^{\Bbb C}((M_{\nabla}^{\Bbb C})_{\iota})=
{\rm dim}\,(F^h)^{-1}({\bf a})$, $\iota^{\Bbb C}((M_{\nabla}^{\Bbb C})_{\iota})$ 
is an open potion of $(F^h)^{-1}({\bf a})$. 
\hspace{10.4truecm}q.e.d.

\vspace{0.5truecm}

\noindent{\it Remark 4.1.} Take another $C^{\omega}$-Koszul connection 
$\widehat{\nabla}$ of $M$.  Let $\hat{\iota}^{\Bbb C}$ be the complexification 
of $\iota$ as a map of $M_{\widehat{\nabla}}^{\Bbb C}$ into 
$\widetilde M_{\widetilde{\nabla}}^{\Bbb C}$.  Take $X\in(M_{\nabla}^{\Bbb C})
_{\iota}\cap(M_{\widehat{\nabla}}^{\Bbb C})_{\iota}\,(\subset TM)$.  
Then $\iota^{\Bbb C}(X)$ and $\hat{\iota^{\Bbb C}}(X)$ are mutually distinct 
in general but they belong to $(F^h)^{-1}({\bf a})$.  

\vspace{0.5truecm}

\noindent
{\it Example.} Let $S^n(r):=\{(x_1,\cdots,x_{n+1})\in{\Bbb R}^{n+1}\,\vert\,
x_1^2+\cdots+x_{n+1}^2=r^2\}$ and $g$ be the standard Riemannian metric of 
$S^n(r)$.  Denote by $\iota$ the inclusion map of $S^n(r)$ into 
${\Bbb R}^{n+1}$.  Then we have 
$$\iota^{\Bbb C}(S^n(r)_g^{\Bbb C})=\{(z_1,\cdots,z_{n+1})\in{\Bbb C}^{n+1}\,
\vert\,z_1^2+\cdots+z_{n+1}^2=r^2\}.$$

\section{The anti-Kaehler metric on the complexification of\\
a pseudo-Riemannian manifold} 
Let $(M,g)$ be an $m$-dimensional $C^{\omega}$-pseudo-Riemannian manifold and 
$M_g^{\Bbb C}=(U^g,J^g)$ be its complexification.  We shall construct an 
anti-Hermitian metric associated with $J^g$ on a tubular neighborhood of $M$ 
in $M_g^{\Bbb C}$.  Fix $p_0\in M$.  Take a holomorphic local coordinate 
$(V,\phi=(z_1,\cdots,z_m))$ of $M_g^{\Bbb C}$ around $p_0$ satisfying 
$M\cap V=\phi^{-1}({\Bbb R}^m)$.  
Let $\phi\vert_{M\cap V}=(x_1,\cdots,x_m)$.  As $g\vert_{M\cap V}=
\sum\limits_{i=1}^n\sum\limits_{j=1}^ng_{ij}dx_idx_j$, we define a holomorphic 
metric $g^{h,p_0}$ on a neighborhood of $M\cap V$ in $V$ by 
$g^{h,p_0}:=\sum\limits_{i=1}^n\sum\limits_{j=1}^ng_{ij}^hdz_idz_j$, where 
$g_{ij}^h$ is a holomorphic extension of $g_{ij}$.  
Thus, for each $p\in M$, we can define a holomorphic metric $g^{h,p}$ on a 
neighborhood of $p$ in $M_g^{\Bbb C}$.  By patching $g^{h,p}$'s ($p\in M$), 
we get a holomorphic metric on a tubular neighborhood of $M$ in 
$M_g^{\Bbb C}$.  Furthermore, we extend this holomorphic metric to the maximal 
one.  Denote by $g^h$ this maximal holomorphic metric.  

\vspace{0.5truecm}

\noindent
{\bf Notation 1.} Denote by $(M_g^{\Bbb C})_A$ the domain of $g^h$.  

\vspace{0.5truecm}

Note that $g^h$ is a holomorphic section of the holomorphic vector bundle 
$(T^{\ast}((M_g^{\Bbb C})_A)\otimes T^{\ast}((M_g^{\Bbb C})_A))^{(2,0)}
(\subset(T^{\ast}((M_g^{\Bbb C})_A)\otimes 
T^{\ast}((M_g^{\Bbb C})_A))^{\Bbb C})$ 
consisting of all complex $(0,2)$-tensors of type $(2,0)$ of 
$(M_g^{\Bbb C})_A$.  
From $g^h$, we define an anti-Kaehler metric associated 
with $J^g$ as follows.  

\vspace{0.3truecm}

\noindent
{\bf Definition 1.} 
Define $\overline{g^h}$ by $\overline{g^h}(Z_1,Z_2)=
\overline{
g^h(\bar Z_1,\bar Z_2)}$ ($Z_1,Z_2\in(T(M_g^{\Bbb C})_A)^{\Bbb C})$, where 
$\overline{(\cdot)}$ is the conjugation of $(\cdot)$.  Then 
$(g^h+\overline{g^h})\vert_{T((M_g^{\Bbb C})_A)\times 
T((M_g^{\Bbb C})_A)}$ is an 
anti-Kaehler metric on $(M_g^{\Bbb C})_A$ (by Theorem 2.2 of [1]).  
We denote this anti-Kaehler metric by $g_A$.  

\vspace{0.5truecm}

\noindent
{\it Remark 5.1.} (i) For $X,Y\in T(M_g^{\Bbb C})$, we have $g_A(X,Y)=2{\rm Re}
(g^h(X,Y))$.  

(ii) If $(M,g)$ is Einstein, then $((M_g^{\Bbb C})_A,g_A)$ also is 
Einstein (see Section 5 of [1]).  Hence 
$((((M_g^{\Bbb C})_A)_{g_A}^{\Bbb C})_A,\,\,(g_A)_A)$ also 
is Einstein.  Thus we get an inductive construction of an Einstein 
(anti-Kaehler) manifold.  

\vspace{0.5truecm}

\noindent
{\bf Notation 2.} For a $C^{\omega}$-map $f:(M,g)\to
(\widetilde M,\widetilde g)$ between $C^{\omega}$-pseudo-Riemannian 
manifolds, we set $(M_g^{\Bbb C})_{A,f:i}:=(M_g^{\Bbb C})_A\cap
(M_g^{\Bbb C})_{f:i}$.  

\vspace{0.5truecm}

For the complexification of a $C^{\omega}$-isometric immersion between 
$C^{\omega}$-pseudo-Riemannain manifolds, we have the following fact.  

\vspace{0.5truecm}

\noindent
{\bf Theorem 5.1.} {\sl Let $f:(M,g)\hookrightarrow
(\widetilde M,\widetilde g)$ be a $C^{\omega}$-isometric immersion between 
$C^{\omega}$-pseudo-Riemannian manifolds.  Then the complexified map 
$f^{\Bbb C}:((M_g^{\Bbb C})_{A,f:i}\cap(f^{\Bbb C})^{-1}
((\widetilde M_{\widetilde g}^{\Bbb C})_A),\,g_A)\to
(\widetilde M_{\widetilde g}^{\Bbb C})_A,\,\widetilde g_A)$ is a holomorphic 
and isometric (that is, anti-Kaehler) immersion.}

\vspace{0.5truecm}

\noindent
{\it Proof.} For simplicity, we set $(M_g^{\Bbb C})'_{A,f:i}:=
(M_g^{\Bbb C})_{A,f:i}\cap(f^{\Bbb C})^{-1}
((\widetilde M_{\widetilde g}^{\Bbb C})_A)$.  
We suffice to show that $(f^{\Bbb C})^{\ast}
\widetilde g_A=g_A$.  Let $g^h$ (resp. $\widetilde g^h$) be a 
holomorphic metric arising from $g$ (resp. $\widetilde g$).  
Since $f^{\Bbb C}$ is holomorphic by Proposition 4.1, 
$((f^{\Bbb C}_{\ast})^{\Bbb C})^{\ast}\widetilde g^h$ is the holomorphic 
$(0,2)$-tensor field on $(M_g^{\Bbb C})'_{A,f:i}$.  
Also, it is clear that $((f^{\Bbb C}_{\ast})^{\Bbb C})^{\ast}\widetilde g^h
\vert_{TM\times TM}=f^{\ast}\widetilde g(=g)$.  Hence we get 
$((f^{\Bbb C}_{\ast})^{\Bbb C})^{\ast}\widetilde g^h=g^h$ on 
$(M_g^{\Bbb C})'_{A,f:i}$ and furthermore 
$$\begin{array}{l}
\displaystyle{(f^{\Bbb C})^{\ast}\widetilde g_A=(f^{\Bbb C})^{\ast}\left(
(\widetilde g^h+\overline{\widetilde g^h})
\vert_{T((\widetilde M_{\widetilde g}^{\Bbb C})_A)\times 
T((\widetilde M^{\Bbb C}_{\widetilde g})_A)}\right)}\\
\hspace{1.6truecm}\displaystyle{=\left(((f^{\Bbb C}_{\ast})^{\Bbb C})^{\ast}
\widetilde g^h+\overline{((f^{\Bbb C}_{\ast})^{\Bbb C})^{\ast}\widetilde g^h}
\right)\vert_{T((M_g^{\Bbb C})'_{A,f:i})\times T((M_g^{\Bbb C})'_{A,f:i})}}\\
\hspace{1.6truecm}\displaystyle{=(g^h+\overline{g^h})
\vert_{T((M_g^{\Bbb C})'_{A,f:i})\times T((M_g^{\Bbb C})'_{A,f:i})}=g_A}
\end{array}$$
on $(M_g^{\Bbb C})'_{A,f:i}$.  
\begin{flushright}q.e.d.\end{flushright}

\vspace{0.5truecm}

\noindent
{\bf Definition 2.} We call the anti-Kaehler submanifold 
$f^{\Bbb C}:((M_g^{\Bbb C})'_{A,f:i},g_A)\hookrightarrow
((\widetilde M_{\widetilde g}^{\Bbb C})_A,\widetilde g_A)$ the 
{\it complexfication of the Riemannian submanifold} 
$f:(M,g)\hookrightarrow(\widetilde M,\widetilde g)$.  

\section{Complete complexifications of pseudo-Riemannian homogeneous spaces}
Let $(G/K,g)$ be a pseudo-Riemannian homogeneous space.  
Here we assume that $G$ and $K$ admit faithful real representations.  
Hence the complexifications $G^{\Bbb C}$ and $K^{\Bbb C}$ of $G$ and $K$ are 
defined.  Since $g_{eK}$ is invariant with respect to the $K$-action on 
$T_{eK}(G/K)$, its complexification $g_{eK}^{\Bbb C}$ is invariant with repsect 
to the $K^{\Bbb C}$-action on $T_{eK^{\Bbb C}}(G^{\Bbb C}/K^{\Bbb C})(=
(T_{eK}(G/K))^{\Bbb C})$.  
Hence we obtain a $G^{\Bbb C}$-invariant holomorphic metric $\widetilde g^h$ 
on $G^{\Bbb C}/K^{\Bbb C}$ from the ${\Bbb C}$-bilinear extension of 
$g_{eK}^{\Bbb C}$ to 
$(T_{eK^{\Bbb C}}(G^{\Bbb C}/K^{\Bbb C}))^{\Bbb C}\times 
(T_{eK^{\Bbb C}}(G^{\Bbb C}/K^{\Bbb C}))^{\Bbb C}$.  
Set $\widetilde g_A:=
2{\rm Re}\widetilde g^h\vert_{T(G^{\Bbb C}/K^{\Bbb C})\times 
T(G^{\Bbb C}/K^{\Bbb C})}$, which is also $G^{\Bbb C}$-inavariant.  Define 
$j:T_{eK^{\Bbb C}}(G^{\Bbb C}/K^{\Bbb C})\to T_{eK^{\Bbb C}}(G^{\Bbb C}/K^{\Bbb C})$ 
by $j(X):=\sqrt{-1}X$ ($X\in T_{eK^{\Bbb C}}(G^{\Bbb C}/K^{\Bbb C})$).  
Since $j$ is invariant with respect to the $K^{\Bbb C}$-action on 
$T_{eK^{\Bbb C}}(G^{\Bbb C}/K^{\Bbb C})$, we obtain a $G^{\Bbb C}$-invariant 
almost complex structure $\widetilde J$ of $G^{\Bbb C}/K^{\Bbb C}$ from $j$.  
Then it is shown that $(\widetilde J,\widetilde g_A)$ is an anti-Kaehler 
structure of $G^{\Bbb C}/K^{\Bbb C}$.  Also, 
it is clear that $(G^{\Bbb C}/K^{\Bbb C},\widetilde J,\widetilde g_A)$ is 
geodesically complete.  By identifying $G/K$ with $G(eK^{\Bbb C})$, 
$G^{\Bbb C}/K^{\Bbb C}$ is regarded as the complete complexification of $G/K$.  
Define $\Phi:T(G/K)\to G^{\Bbb C}/K^{\Bbb C}$ by $\Phi(v):=\exp_p
(\widetilde J_pv)$ for $v\in T(G/K)$, where $p$ is the base point of $v$ and 
$\exp_p$ is the exponential map of the anti-Kaehler manifold 
$(G^{\Bbb C}/K^{\Bbb C},\widetilde J,\widetilde g_A)$ at $p\,(\in G/K
=G(eK^{\Bbb C})\subset G^{\Bbb C}/K^{\Bbb C})$.  
Note that this map $\Phi$ is called the polar map in [5].  

\vspace{0.5truecm}

\noindent
{\it Remark 6.1.} For a $C^{\omega}$-isometric immersion $f$ of a 
$C^{\omega}$-Riemannian manifold $(M,g)$ into a Riemannian symmetric space 
$(G/K,g)$ of non-compact type, we [15] defined its complexification 
as an immersion of a tubular neighborhood of $M$ in $(M_g^{\Bbb C})_{f:i}$ into 
$G^{\Bbb C}/K^{\Bbb C}$.  It is shown that the complexification defined in [15] 
is equal to the composition of the complexification 
$f^{\Bbb C}(:(M_g^{\Bbb C})_{f:i}\to(G/K)_g^{\Bbb C})$ defined in 
Section 4 and the polar map $\Phi$.  

\vspace{0.5truecm}

Set $\widetilde{\Omega}:=\displaystyle{
\mathop{\cup}_{v\in T^{\perp}G(eK^{\Bbb C})}\{\exp(sv)\,\vert\,
0\leq s<r_v\}}$, where $\exp$ is the exponential map of 
$G^{\Bbb C}/K^{\Bbb C}$ and $r_v$ is the first focal radius of 
$G(eK^{\Bbb C})(\subset G^{\Bbb C}/K^{\Bbb C})$ 
along $\gamma_v$.  
We have the following fact for $\Phi$.  

\vspace{0.5truecm}

\noindent
{\bf Theorem 6.1.} {\sl The restriction $\Phi\vert_{((G/K)_g^{\Bbb C})_A}$ of 
$\Phi$ to $((G/K)_g^{\Bbb C})_A$ is a diffeomorphism onto $\widetilde{\Omega}$ 
and, each point of the boundary $\partial((G/K)_g^{\Bbb C})_A$ of 
$((G/K)_g^{\Bbb C})_A$ in $T(G/K)$ is a critical point of $\Phi$.  Furthermore, 
$\Phi\vert_{((G/K)_g^{\Bbb C})_A}$ is a holomorphic isometry (that is, 
$(\Phi\vert_{((G/K)_g^{\Bbb C})_A})^{\ast}\widetilde J=J^g$ and 
$(\Phi\vert_{((G/K)_g^{\Bbb C})_A})^{\ast}\widetilde g_A=g_A$).}

\vspace{0.5truecm}

\noindent
{\it Proof.} Let $\Omega$ be the connected component of $T(G/K)$ containing 
the $0$-section $(=G/K)$ of the set of all regular points of $\Phi$.  From 
the definition of $\Phi$, it is easy to show that $v\in T_p(G/K)(\subset 
T(G/K))$ is a critical point of $\Phi$ if and only if $\Phi(v)$ is a focal 
point of the orbit $G(eK^{\Bbb C})$ along $\gamma_v$ or a conjugate point of 
$p$ along $\gamma_v$.  Hence we see that $\Phi(\Omega)=\widetilde{\Omega}$ and 
that $\Phi\vert_{\Omega}$ is a diffeomorphism onto $\widetilde{\Omega}$.  
Now we shall show that $\Phi\vert_{\Omega}$ is a holomorphic isometry.  Let 
$\gamma$ be a geodesic in $G/K$.  We have 
$$\begin{array}{l}
\displaystyle{(\Phi\circ\gamma_{\ast})(s+t\sqrt{-1})=\Phi(t\gamma'(s))
=\exp_{\gamma(s)}(\widetilde J_{\gamma(s)}(t\gamma'(s)))}\\
\hspace{3.4truecm}\displaystyle{=(\gamma_{t\gamma'(s)})^{\Bbb C}(\sqrt{-1})
=\gamma^{\Bbb C}(s+t\sqrt{-1}),}
\end{array}$$
where $(\gamma_{t\gamma'(s)})^{\Bbb C}$ (resp. $\gamma^{\Bbb C}$) is 
the complexification of $\gamma_{t\gamma'(s)}$ (resp. $\gamma$) in 
$G^{\Bbb C}/K^{\Bbb C}$.  Thus $\Phi\circ\gamma_{\ast}(:T{\Bbb R}={\Bbb C}\to
(G^{\Bbb C}/K^{\Bbb C},\widetilde J))$ is holomorphic.  Therefore, according to 
Theorem 3.4 of [27], $\Phi\vert_{(G/K)_g^{\Bbb C}}$ is holomorphic, that is, 
$(\Phi\vert_{(G/K)_g^{\Bbb C}})^{\ast}\widetilde J=J_A$.  
On the other hand, it is clear that $(\Phi\vert_{\Omega})^{\ast}\widetilde J$ 
is equal to $J_A$ on $\Omega$.  Hence we have $\Omega\subset(G/K)_g^{\Bbb C}$.  
Since $(\Phi\vert_{\Omega})^{\ast}\widetilde g^h$ is 
the non-extendable holomorphic metric arising from $g$.  Hence we have 
$\Omega=((G/K)_g^{\Bbb C})_A$.  Hence 
the statement of this theorem follows.  \hspace{8truecm} q.e.d.

\vspace{0.5truecm}

\section{Duals of a pseudo-Riemannian manifolds}
In this section, we shall define the dual of a 
$C^{\omega}$-pseudo-Riemannian manifold and the dual of a totally 
geodesic $C^{\omega}$-pseudo-Riemannian submanifold.  
Let $(M,g)$ be a $C^{\omega}$-pseudo-Riemannian manifold.  
For each $p\in M$, we set $M_p^{\ast}:=
(M_g^{\Bbb C})_A\cap T_pM$ and denote the inclusion map of $M^{\ast}_p$ 
into $(M_g^{\Bbb C})_A$ by $\iota_p$.  For $M_p^{\ast}$, the following fact 
holds.  



\centerline{
\unitlength 0.1in
\begin{picture}( 60.3000, 38.9000)( -2.1000,-40.6000)
%
\special{pn 8}%
\special{pa 2000 2000}%
\special{pa 3410 2000}%
\special{fp}%
%
\special{pn 8}%
\special{pa 2000 1400}%
\special{pa 3410 1400}%
\special{pa 3410 2600}%
\special{pa 2000 2600}%
\special{pa 2000 1400}%
\special{dt 0.045}%
%
\special{pn 8}%
\special{pa 2000 1400}%
\special{pa 3410 1400}%
\special{pa 3410 2000}%
\special{pa 2000 2000}%
\special{pa 2000 1400}%
\special{dt 0.045}%
%
\special{pn 8}%
\special{pa 2710 1200}%
\special{pa 2710 2800}%
\special{fp}%
\put(45.0000,-7.9000){\makebox(0,0)[rt]{$0$-section$=G/K$}}%
\put(30.3000,-40.5000){\makebox(0,0)[rt]{$T(G/K)$}}%
%
\special{pn 8}%
\special{ar 1950 2000 230 600  1.5707963 1.5997120}%
\special{ar 1950 2000 230 600  1.6864590 1.7153746}%
\special{ar 1950 2000 230 600  1.8021216 1.8310373}%
\special{ar 1950 2000 230 600  1.9177843 1.9466999}%
\special{ar 1950 2000 230 600  2.0334469 2.0623626}%
\special{ar 1950 2000 230 600  2.1491096 2.1780252}%
\special{ar 1950 2000 230 600  2.2647722 2.2936879}%
\special{ar 1950 2000 230 600  2.3804349 2.4093505}%
\special{ar 1950 2000 230 600  2.4960975 2.5250132}%
\special{ar 1950 2000 230 600  2.6117602 2.6406758}%
\special{ar 1950 2000 230 600  2.7274228 2.7563385}%
\special{ar 1950 2000 230 600  2.8430855 2.8720011}%
\special{ar 1950 2000 230 600  2.9587481 2.9876638}%
\special{ar 1950 2000 230 600  3.0744108 3.1033264}%
\special{ar 1950 2000 230 600  3.1900734 3.2189891}%
\special{ar 1950 2000 230 600  3.3057361 3.3346517}%
\special{ar 1950 2000 230 600  3.4213987 3.4503144}%
\special{ar 1950 2000 230 600  3.5370614 3.5659770}%
\special{ar 1950 2000 230 600  3.6527240 3.6816397}%
\special{ar 1950 2000 230 600  3.7683867 3.7973024}%
\special{ar 1950 2000 230 600  3.8840493 3.9129650}%
\special{ar 1950 2000 230 600  3.9997120 4.0286277}%
\special{ar 1950 2000 230 600  4.1153746 4.1442903}%
\special{ar 1950 2000 230 600  4.2310373 4.2599530}%
\special{ar 1950 2000 230 600  4.3466999 4.3756156}%
\special{ar 1950 2000 230 600  4.4623626 4.4912783}%
\special{ar 1950 2000 230 600  4.5780252 4.6069409}%
\special{ar 1950 2000 230 600  4.6936879 4.7123890}%
\put(16.8000,-18.0000){\makebox(0,0)[rt]{$((G/K)_g^{\Bbb C})_A$}}%
%
\special{pn 4}%
\special{pa 2000 1400}%
\special{pa 3410 1400}%
\special{pa 3410 2600}%
\special{pa 2000 2600}%
\special{pa 2000 1400}%
\special{ip}%
%
\special{pn 8}%
\special{pa 2910 1200}%
\special{pa 2910 2800}%
\special{fp}%
%
\special{pn 8}%
\special{pa 3110 1200}%
\special{pa 3110 2800}%
\special{fp}%
%
\special{pn 8}%
\special{pa 2510 1200}%
\special{pa 2510 2800}%
\special{fp}%
%
\special{pn 8}%
\special{pa 2310 1200}%
\special{pa 2310 2800}%
\special{fp}%
%
\special{pn 8}%
\special{pa 2310 2000}%
\special{pa 2380 2000}%
\special{pa 2380 1940}%
\special{pa 2310 1940}%
\special{pa 2310 2000}%
\special{fp}%
%
\special{pn 8}%
\special{pa 2510 2000}%
\special{pa 2580 2000}%
\special{pa 2580 1940}%
\special{pa 2510 1940}%
\special{pa 2510 2000}%
\special{fp}%
%
\special{pn 8}%
\special{pa 2710 2000}%
\special{pa 2780 2000}%
\special{pa 2780 1940}%
\special{pa 2710 1940}%
\special{pa 2710 2000}%
\special{fp}%
%
\special{pn 8}%
\special{pa 2910 2000}%
\special{pa 2980 2000}%
\special{pa 2980 1940}%
\special{pa 2910 1940}%
\special{pa 2910 2000}%
\special{fp}%
%
\special{pn 8}%
\special{pa 3110 2000}%
\special{pa 3180 2000}%
\special{pa 3180 1940}%
\special{pa 3110 1940}%
\special{pa 3110 2000}%
\special{fp}%
\put(39.9000,-17.5000){\makebox(0,0)[rt]{$\Phi$}}%
\put(40.6000,-19.3000){\makebox(0,0)[rt]{$\longrightarrow$}}%
%
\special{pn 8}%
\special{pa 4400 2000}%
\special{pa 5810 2000}%
\special{fp}%
%
\special{pn 4}%
\special{pa 4400 1400}%
\special{pa 5800 1400}%
\special{dt 0.027}%
%
\special{pn 4}%
\special{pa 4400 2600}%
\special{pa 5800 2600}%
\special{dt 0.027}%
%
\special{pn 8}%
\special{pa 5100 3310}%
\special{pa 5100 710}%
\special{fp}%
%
\special{pn 8}%
\special{ar 5000 2000 290 640  5.0726923 6.2831853}%
\special{ar 5000 2000 290 640  0.0000000 1.1856388}%
%
\special{pn 8}%
\special{ar 5000 810 290 640  5.0726923 6.2831853}%
\special{ar 5000 810 290 640  0.0000000 1.1856388}%
%
\special{pn 8}%
\special{ar 5000 3200 290 640  5.0726923 6.2831853}%
\special{ar 5000 3200 290 640  0.0000000 1.1856388}%
%
\special{pn 8}%
\special{ar 5200 2000 290 640  1.9559538 4.3520857}%
%
\special{pn 8}%
\special{ar 5200 810 290 640  1.9559538 4.3520857}%
%
\special{pn 8}%
\special{ar 5200 3200 290 640  1.9559538 4.3520857}%
%
\special{pn 8}%
\special{pa 5100 790}%
\special{pa 5100 170}%
\special{fp}%
%
\special{pn 8}%
\special{pa 5100 3270}%
\special{pa 5100 3870}%
\special{fp}%
%
\special{pn 8}%
\special{ar 4750 2000 750 670  5.1984668 6.2831853}%
\special{ar 4750 2000 750 670  0.0000000 1.0938903}%
%
\special{pn 8}%
\special{ar 4760 800 750 670  5.1984668 6.2831853}%
\special{ar 4760 800 750 670  0.0000000 1.0938903}%
%
\special{pn 8}%
\special{ar 4760 3190 750 670  5.1984668 6.2831853}%
\special{ar 4760 3190 750 670  0.0000000 1.0938903}%
%
\special{pn 8}%
\special{ar 5440 2000 750 670  2.0477024 4.2263112}%
%
\special{pn 8}%
\special{ar 5440 800 750 670  2.0477024 4.2263112}%
%
\special{pn 8}%
\special{ar 5440 3190 750 670  2.0477024 4.2263112}%
%
\special{pn 8}%
\special{pa 5100 2000}%
\special{pa 5170 2000}%
\special{pa 5170 1940}%
\special{pa 5100 1940}%
\special{pa 5100 2000}%
\special{fp}%
%
\special{pn 8}%
\special{pa 5290 2000}%
\special{pa 5360 2000}%
\special{pa 5360 1940}%
\special{pa 5290 1940}%
\special{pa 5290 2000}%
\special{fp}%
%
\special{pn 8}%
\special{pa 5500 2000}%
\special{pa 5570 2000}%
\special{pa 5570 1940}%
\special{pa 5500 1940}%
\special{pa 5500 2000}%
\special{fp}%
%
\special{pn 8}%
\special{pa 4920 2000}%
\special{pa 4990 2000}%
\special{pa 4990 1940}%
\special{pa 4920 1940}%
\special{pa 4920 2000}%
\special{fp}%
%
\special{pn 8}%
\special{pa 4700 2000}%
\special{pa 4770 2000}%
\special{pa 4770 1940}%
\special{pa 4700 1940}%
\special{pa 4700 2000}%
\special{fp}%
%
\special{pn 4}%
\special{pa 1990 1400}%
\special{pa 3400 1400}%
\special{pa 3400 2600}%
\special{pa 1990 2600}%
\special{pa 1990 1400}%
\special{ip}%
\put(53.1000,-40.6000){\makebox(0,0)[rt]{$G^{\Bbb C}/K^{\Bbb C}$}}%
\put(48.3000,-2.0000){\makebox(0,0)[rt]{$G(eK^{\Bbb C})(=G/K)$}}%
%
\special{pn 8}%
\special{pa 2000 200}%
\special{pa 3410 200}%
\special{pa 3410 3800}%
\special{pa 2000 3800}%
\special{pa 2000 200}%
\special{fp}%
%
\special{pn 8}%
\special{pa 2310 2770}%
\special{pa 2310 3800}%
\special{fp}%
%
\special{pn 8}%
\special{pa 2310 1270}%
\special{pa 2310 200}%
\special{fp}%
%
\special{pn 8}%
\special{pa 2510 1250}%
\special{pa 2510 200}%
\special{fp}%
%
\special{pn 8}%
\special{pa 2710 1270}%
\special{pa 2710 200}%
\special{fp}%
%
\special{pn 8}%
\special{pa 2910 1290}%
\special{pa 2910 200}%
\special{fp}%
%
\special{pn 8}%
\special{pa 3110 1240}%
\special{pa 3110 200}%
\special{fp}%
%
\special{pn 8}%
\special{pa 2510 2730}%
\special{pa 2510 3800}%
\special{fp}%
%
\special{pn 8}%
\special{pa 2710 2720}%
\special{pa 2710 3800}%
\special{fp}%
%
\special{pn 8}%
\special{pa 2910 2770}%
\special{pa 2910 3800}%
\special{fp}%
%
\special{pn 8}%
\special{pa 3110 2730}%
\special{pa 3110 3800}%
\special{fp}%
%
\special{pn 8}%
\special{pa 4470 600}%
\special{pa 4630 2000}%
\special{dt 0.045}%
\special{sh 1}%
\special{pa 4630 2000}%
\special{pa 4642 1932}%
\special{pa 4624 1948}%
\special{pa 4604 1936}%
\special{pa 4630 2000}%
\special{fp}%
%
\special{pn 8}%
\special{pa 3790 980}%
\special{pa 3260 2000}%
\special{dt 0.045}%
\special{sh 1}%
\special{pa 3260 2000}%
\special{pa 3308 1950}%
\special{pa 3286 1954}%
\special{pa 3274 1932}%
\special{pa 3260 2000}%
\special{fp}%
%
\special{pn 8}%
\special{pa 3990 400}%
\special{pa 4460 580}%
\special{dt 0.045}%
%
\special{pn 4}%
\special{pa 1990 1390}%
\special{pa 1990 2600}%
\special{pa 3400 2600}%
\special{pa 3400 1390}%
\special{pa 3400 1390}%
\special{pa 1990 1390}%
\special{ip}%
%
\special{pn 4}%
\special{pa 1990 1450}%
\special{pa 2050 1390}%
\special{fp}%
\special{pa 1990 1510}%
\special{pa 2110 1390}%
\special{fp}%
\special{pa 1990 1570}%
\special{pa 2170 1390}%
\special{fp}%
\special{pa 1990 1630}%
\special{pa 2230 1390}%
\special{fp}%
\special{pa 1990 1690}%
\special{pa 2290 1390}%
\special{fp}%
\special{pa 1990 1750}%
\special{pa 2350 1390}%
\special{fp}%
\special{pa 1990 1810}%
\special{pa 2410 1390}%
\special{fp}%
\special{pa 1990 1870}%
\special{pa 2470 1390}%
\special{fp}%
\special{pa 1990 1930}%
\special{pa 2530 1390}%
\special{fp}%
\special{pa 1990 1990}%
\special{pa 2590 1390}%
\special{fp}%
\special{pa 1990 2050}%
\special{pa 2650 1390}%
\special{fp}%
\special{pa 1990 2110}%
\special{pa 2710 1390}%
\special{fp}%
\special{pa 1990 2170}%
\special{pa 2770 1390}%
\special{fp}%
\special{pa 1990 2230}%
\special{pa 2830 1390}%
\special{fp}%
\special{pa 1990 2290}%
\special{pa 2890 1390}%
\special{fp}%
\special{pa 1990 2350}%
\special{pa 2950 1390}%
\special{fp}%
\special{pa 1990 2410}%
\special{pa 3010 1390}%
\special{fp}%
\special{pa 1990 2470}%
\special{pa 3070 1390}%
\special{fp}%
\special{pa 1990 2530}%
\special{pa 3130 1390}%
\special{fp}%
\special{pa 1990 2590}%
\special{pa 3190 1390}%
\special{fp}%
\special{pa 2040 2600}%
\special{pa 3250 1390}%
\special{fp}%
\special{pa 2100 2600}%
\special{pa 3310 1390}%
\special{fp}%
\special{pa 2160 2600}%
\special{pa 3370 1390}%
\special{fp}%
\special{pa 2220 2600}%
\special{pa 3400 1420}%
\special{fp}%
\special{pa 2280 2600}%
\special{pa 3400 1480}%
\special{fp}%
\special{pa 2340 2600}%
\special{pa 3400 1540}%
\special{fp}%
\special{pa 2400 2600}%
\special{pa 3400 1600}%
\special{fp}%
\special{pa 2460 2600}%
\special{pa 3400 1660}%
\special{fp}%
\special{pa 2520 2600}%
\special{pa 3400 1720}%
\special{fp}%
\special{pa 2580 2600}%
\special{pa 3400 1780}%
\special{fp}%
%
\special{pn 4}%
\special{pa 2640 2600}%
\special{pa 3400 1840}%
\special{fp}%
\special{pa 2700 2600}%
\special{pa 3400 1900}%
\special{fp}%
\special{pa 2760 2600}%
\special{pa 3400 1960}%
\special{fp}%
\special{pa 2820 2600}%
\special{pa 3400 2020}%
\special{fp}%
\special{pa 2880 2600}%
\special{pa 3400 2080}%
\special{fp}%
\special{pa 2940 2600}%
\special{pa 3400 2140}%
\special{fp}%
\special{pa 3000 2600}%
\special{pa 3400 2200}%
\special{fp}%
\special{pa 3060 2600}%
\special{pa 3400 2260}%
\special{fp}%
\special{pa 3120 2600}%
\special{pa 3400 2320}%
\special{fp}%
\special{pa 3180 2600}%
\special{pa 3400 2380}%
\special{fp}%
\special{pa 3240 2600}%
\special{pa 3400 2440}%
\special{fp}%
\special{pa 3300 2600}%
\special{pa 3400 2500}%
\special{fp}%
\special{pa 3360 2600}%
\special{pa 3400 2560}%
\special{fp}%
%
\special{pn 4}%
\special{pa 4410 1430}%
\special{pa 4440 1400}%
\special{fp}%
\special{pa 4410 1490}%
\special{pa 4500 1400}%
\special{fp}%
\special{pa 4410 1550}%
\special{pa 4560 1400}%
\special{fp}%
\special{pa 4410 1610}%
\special{pa 4620 1400}%
\special{fp}%
\special{pa 4410 1670}%
\special{pa 4680 1400}%
\special{fp}%
\special{pa 4410 1730}%
\special{pa 4740 1400}%
\special{fp}%
\special{pa 4410 1790}%
\special{pa 4800 1400}%
\special{fp}%
\special{pa 4410 1850}%
\special{pa 4860 1400}%
\special{fp}%
\special{pa 4410 1910}%
\special{pa 4920 1400}%
\special{fp}%
\special{pa 4410 1970}%
\special{pa 4980 1400}%
\special{fp}%
\special{pa 4410 2030}%
\special{pa 5040 1400}%
\special{fp}%
\special{pa 4410 2090}%
\special{pa 5100 1400}%
\special{fp}%
\special{pa 4410 2150}%
\special{pa 5160 1400}%
\special{fp}%
\special{pa 4410 2210}%
\special{pa 5220 1400}%
\special{fp}%
\special{pa 4410 2270}%
\special{pa 5280 1400}%
\special{fp}%
\special{pa 4410 2330}%
\special{pa 5340 1400}%
\special{fp}%
\special{pa 4410 2390}%
\special{pa 5400 1400}%
\special{fp}%
\special{pa 4410 2450}%
\special{pa 5460 1400}%
\special{fp}%
\special{pa 4410 2510}%
\special{pa 5520 1400}%
\special{fp}%
\special{pa 4410 2570}%
\special{pa 5580 1400}%
\special{fp}%
\special{pa 4440 2600}%
\special{pa 5640 1400}%
\special{fp}%
\special{pa 4500 2600}%
\special{pa 5700 1400}%
\special{fp}%
\special{pa 4560 2600}%
\special{pa 5760 1400}%
\special{fp}%
\special{pa 4620 2600}%
\special{pa 5820 1400}%
\special{fp}%
\special{pa 4680 2600}%
\special{pa 5820 1460}%
\special{fp}%
\special{pa 4740 2600}%
\special{pa 5820 1520}%
\special{fp}%
\special{pa 4800 2600}%
\special{pa 5820 1580}%
\special{fp}%
\special{pa 4860 2600}%
\special{pa 5820 1640}%
\special{fp}%
\special{pa 4920 2600}%
\special{pa 5820 1700}%
\special{fp}%
\special{pa 4980 2600}%
\special{pa 5820 1760}%
\special{fp}%
%
\special{pn 4}%
\special{pa 5040 2600}%
\special{pa 5820 1820}%
\special{fp}%
\special{pa 5100 2600}%
\special{pa 5820 1880}%
\special{fp}%
\special{pa 5160 2600}%
\special{pa 5820 1940}%
\special{fp}%
\special{pa 5220 2600}%
\special{pa 5820 2000}%
\special{fp}%
\special{pa 5280 2600}%
\special{pa 5820 2060}%
\special{fp}%
\special{pa 5340 2600}%
\special{pa 5820 2120}%
\special{fp}%
\special{pa 5400 2600}%
\special{pa 5820 2180}%
\special{fp}%
\special{pa 5460 2600}%
\special{pa 5820 2240}%
\special{fp}%
\special{pa 5520 2600}%
\special{pa 5820 2300}%
\special{fp}%
\special{pa 5580 2600}%
\special{pa 5820 2360}%
\special{fp}%
\special{pa 5640 2600}%
\special{pa 5820 2420}%
\special{fp}%
\special{pa 5700 2600}%
\special{pa 5820 2480}%
\special{fp}%
\special{pa 5760 2600}%
\special{pa 5820 2540}%
\special{fp}%
\end{picture}%
\hspace{3.2truecm}}

\vspace{0.5truecm}

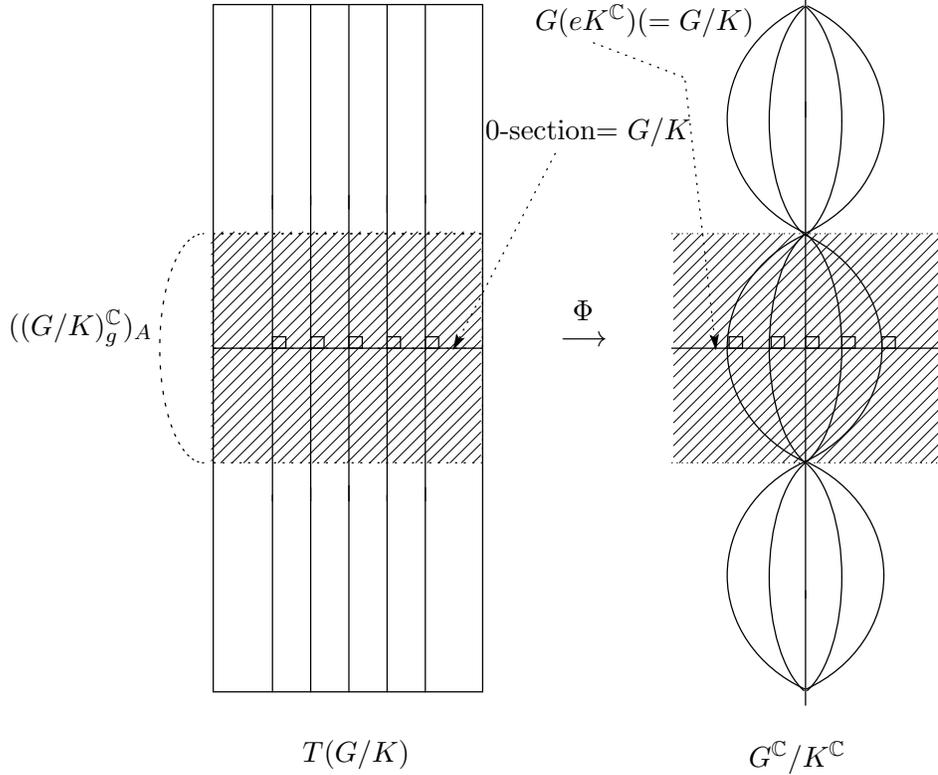
\centerline{{\bf Figure 2.}}

\vspace{0.5truecm}

\noindent
{\bf Proposition 7.1.} {\sl Let $\exp_p$ be the exponential map of 
$((M_g^{\Bbb C})_A,g_A)$ at $p$ and $D_p\,(\subset T_p((M_g^{\Bbb C})_A))$ be 
its domain.  The above set $M_p^{\ast}$ coincides with the geodesic umbrella 
$\exp_p(T_p(M_p^{\ast})\cap D)$.}

\vspace{0.5truecm}

\noindent
{\it Proof.} For each $X\in M_p^{\ast}$, we get 
${\rm id}_M^{\Bbb C}(X)=\gamma_X^{\Bbb C}(\sqrt{-1})=\exp_p(J^g_pX)$.  
On the other hand, it is clear that ${\rm id}_M^{\Bbb C}
={\rm id}_{M_g^{\Bbb C}}$.  Hence we get $X=\exp_p(J^g_pX)\in
\exp_p(T_p(M_p^{\ast})\cap D)$.  From the arbitrariness of $X$, we get 
$M_p^{\ast}\subset\exp_p(T_p(M_p^{\ast})\cap D)$.  
It is clear that this relation implies 
$M_p^{\ast}=\exp_p(T_p(M_p^{\ast})\cap D)$.
\hspace{9.7truecm}q.e.d.

\vspace{0.5truecm}

\noindent
{\bf Definition 3.} We call the pseudo-Riemannian manifold $(M_p^{\ast},
\iota_p^{\ast}g_A)$ the {\it dual of} $(M,g)$ {\it at} $p$.  

\vspace{0.5truecm}

The following question is proposed naturally:  

\vspace{0.5truecm}

{\sl Are $(M,g)$ and $(M_p^{\ast},\iota_p^{\ast}g_A)$ 
totally geodesic in $((M_g^{\Bbb C})_A,g_A)$?}

\vspace{0.5truecm}

For this question, we can show the following fact.  

\vspace{0.5truecm}

\noindent
{\bf Proposition 7.2.} {\sl The submanifold $(M,g)$ is totally geodesic in 
$((M_g^{\Bbb C})_A,g_A)$.}

\vspace{0.5truecm}

\noindent
{\it Proof.} Define $\sigma:M_g^{\Bbb C}\to M_g^{\Bbb C}$ by $\sigma(X)=-X$ 
($X\in (M_g^{\Bbb C})_A$).  It is clear that $\sigma$ is an isometry of 
$((M_g^{\Bbb C})_A,g_A)$.  Hence, since $M$ is a component of the fixed point 
set of $\sigma$, $(M,g)$ is totally geodesic in $((M_g^{\Bbb C})_A,g_A)$.  
\hspace{9.8truecm}q.e.d.

\vspace{0.5truecm}

Also, we can show the following fact in the case where $(M,g)$ is a 
pseudo-Riemannian symmetric space.  

\vspace{0.5truecm}

\noindent
{\bf Theorem 7.3.} {\sl Let $(G/K,g)$ be a pseudo-Riemannian 
symmetric space associated with a semi-simple symmetric pair $(G,K)$.  
Then $((G/K)_p^{\ast},\iota_p^{\ast}g_A)$ is totally geodesic in \newline
$(((G/K)_g^{\Bbb C})_A,g_A)$.}

\vspace{0.5truecm}

\noindent
{\it Proof.} We suffice to show the statement in case of $p=eK(=eK^{\Bbb C})$ 
($e\,:\,$ the identity element of $G$).  Let $\mathfrak g$ be the Lie algebra 
of $G$ and $\mathfrak g=\mathfrak k+\mathfrak p$ be the canonical 
decomposition associated with $(G,K)$.  Then $T_{eK}(G^{\Bbb C}/K^{\Bbb C})$ 
is identified with $\mathfrak p^{\Bbb C}$.  
Let $\Phi$ be as in Section 6.  It follows from the definition of $\Phi$ that 
$\exp_{eK^{\Bbb C}}(\sqrt{-1}\mathfrak p)\supset\Phi((G/K)^{\ast}_{eK})$.  
Since $\sqrt{-1}\mathfrak p$ is a Lie triple system of 
$\mathfrak p^{\Bbb C}$, $\exp_{eK^{\Bbb C}}(\sqrt{-1}\mathfrak p)$ 
is totally geodesic in $G^{\Bbb C}/K^{\Bbb C}$.  
Hence, since $\Phi\vert_{((G/K)_g^{\Bbb C})_A}$ is an isometry into 
$G^{\Bbb C}/K^{\Bbb C}$ by Theorem 6.1, $(G/K)^{\ast}_{eK}$ is totally 
geodesic in $(((G/K)_g^{\Bbb C})_A,g_A)$.  \hspace{11.1truecm}q.e.d.

\vspace{0.5truecm}

Let $f:(M,g)\hookrightarrow(\widetilde M,\widetilde g)$ be a 
$C^{\omega}$-isometric immersion between $C^{\omega}$-pseudo-Riemannian 
manifolds and set $(M_p^{\ast})_f:=M_p^{\ast}\cap(M_g^{\Bbb C})_f$.  
Then the following question is proposed naturally:

\vspace{0.2truecm}

{\sl Is $f^{\Bbb C}((M_p^{\ast})_f)$ contained in 
${\widetilde M}^{\ast}_{f(p)}$ for each $p\in M$?}

\vspace{0.2truecm}

\noindent
For this problem, we have the following fact.  

\vspace{0.5truecm}

\noindent
{\bf Theorem 7.4.} {\sl If $f$ is totally geodesic, then $f^{\Bbb C}
((M_p^{\ast})_f)$ is contained in ${\widetilde M}^{\ast}_{f(p)}$ for each 
$p\in M$.}

\vspace{0.5truecm}

\noindent
{\it Proof.} Let $X\in(M_p^{\ast})_f$.  Denote by $\exp_{f(p)}$ 
the exponential map of $((\widetilde M_{\widetilde g}^{\Bbb C})_A,
\widetilde g_A)$ at $f(p)$.  Since $f$ is totally geodesic and 
$\exp_{f(p)}$ is holomorphic, we have 
$$\begin{array}{r}
\displaystyle{f^{\Bbb C}(X)=(f\circ\gamma_X)^h(\sqrt{-1})
=(\gamma_{f_{\ast}(X)})^{\Bbb C}(\sqrt{-1})
=\exp_{f(p)}(J^{\widetilde g}_{f(p)}(f_{\ast}(X)))}\\
\hspace{3.5truecm}\displaystyle{\in\exp_{f(p)}(T_{f(p)}
\widetilde M^{\ast}_{f(p)}\cap D),}
\end{array}$$
where $\gamma_X$ (resp. $\gamma_{f_{\ast}(X)}$) is the geodesic in $(M,g)$ 
(resp. $(\widetilde M,\widetilde g)$) with $\dot{\gamma}_X(0)=X$ (resp. 
$\dot{\gamma}_{f_{\ast}(X)}(0)=f_{\ast}(X)$) and $D$ is the domain of 
$\exp_{f(p)}$.  According to Proposition 7.1, 
$\exp_{f(p)}(T_{f(p)}\widetilde M_{f(p)}^{\ast}\cap D)$ is equal to 
${\widetilde M}^{\ast}_{f(p)}$.  Therefore, 
we get $f^{\Bbb C}((M_p^{\ast})_f)\subset{\widetilde M}^{\ast}_{f(p)}$.  
\begin{flushright}q.e.d.\end{flushright}

\vspace{0.5truecm}

\noindent
{\bf Definition 4.} For a totally geodesic $C^{\omega}$-pseudo-Riemannian 
submanifold $f(M)$ in $(\widetilde M,\widetilde g)$, we call a submanifold 
$f^{\Bbb C}((M^{\ast}_p)_f)$ in $({\widetilde M}^{\ast}_{f(p)},
\iota^{\ast}_{f(p)}\widetilde g_A)$ the {\it dual of} $f(M)$.  

\vspace{0.5truecm}

\noindent
{\it Example.} Let $G/K$ be a pseudo-Riemannian symmetric space, $H$ be 
a symmetric subgroup of $G$, $\theta$ be the involution of $G$ with 
$({\rm Fix}\,\theta)_0\subset K\subset {\rm Fix}\,\theta$ and $\sigma$ be 
the involution of $G$ with $({\rm Fix}\,\sigma)_0\subset H\subset {\rm Fix}\,
\sigma$, where $({\rm Fix}\,\theta)_0$ (resp. $({\rm Fix}\,\sigma)_0$) is 
the identity component of ${\rm Fix}\,\theta$ (resp. ${\rm Fix}\,\sigma$).  
Assume that $\theta\circ\sigma=\sigma\circ\theta$.  Also, let $G^{\ast}$ be 
the dual of $G$ with respect to $K$ and $H^{\ast}$ be the dual of $H$ with 
respect to $H\cap K$.  Then the orbit $H(eK)\,(\subset G/K)$ is totally 
geodesic and hence $\iota^{\Bbb C}((H(eK))_{eK}^{\ast})$ is contained in 
$(G/K)_{eK}^{\ast}(=G^{\ast}/K)$, where $\iota^{\Bbb C}$ is the 
complexification of the inclusion map of $H(eK)$ into $G/K$.  Furthermore, 
$\iota^{\Bbb C}((H(eK))_{eK}^{\ast})$ coincides with the orbit 
$H^{\ast}(eK)\,(\subset G^{\ast}/K=(G/K)^{\ast}_{eK}$).   

\section{Complex focal radii} 
In this section, we shall introduce the notions of a complex Jacobi field 
along a complex geodesic in an anti-Kaehler manifold.  
Also, we give a new definition of a complex focal radius of anti-Kaehler 
submanifold by using the notion of a complex Jacobi field and show that 
the notion of a complex focal radius by this new definition coincides with 
one defined in [15] (see Proposition 8.4).  
Next we show a fact which is very useful to calculate the complex focal 
radii of an anti-Kaehler submanifold with section in an anti-Kaehler 
symmetric space (see Proposition 8.5).  
Also, we show that a complex focal radius of a $C^{\omega}$-Riemannian 
submanifold in a Riemannian symmetric space $G/K$ of non-compact type 
(see Definition 6 about the definition of this notion) coincides with 
one defined in [14] (see Proposition 8.6).  
Let $(M,J,g)$ be an anti-Kaehler manifold, $\nabla$ (resp. $R$) be the 
Levi-Civita connection (resp. the curvature tensor) of $g$ 
and $\nabla^{\Bbb C}$ (resp. $R^{\Bbb C}$) be the complexification of 
$\nabla$ (resp. $R$).  Let $(TM)^{(1,0)}$ be the holomorphic vector bundle 
consisting of complex vectors of $M$ of type $(1,0)$.  Note that the 
restriction of $\nabla^{\Bbb C}$ to $TM^{(1,0)}$ is a holomorphic connection 
of $TM^{(1,0)}$ (see Theorem 2.2 of [1]).  
For simplicity, assume that $(M,J,g)$ is complete even if 
the discussion of this section is valid without the assumption of 
the completeness of $(M,J,g)$.  
Let $\gamma:{\Bbb C}\to M$ be 
a complex geodesic, that is, $\gamma(z)=\exp_{\gamma(0)}(({\rm Re}\,z)
\gamma_{\ast}((\frac{\partial}{\partial s})_0)+({\rm Im}\,z)
J_{\gamma(0)}\gamma_{\ast}((\frac{\partial}{\partial s})_0))$, where $(z)$ is 
the complex coordinate of ${\Bbb C}$ and $s:={\rm Re}\,z$.  
Let $Y:{\Bbb C}\to (TM)^{(1,0)}$ be a holomorphic vector 
field along $\gamma$.  
That is, $Y$ assigns $Y_z\in(T_{\gamma(z)}M)^{(1,0)}$ to each 
$z\in {\Bbb C}$ and, for each holomorphic local coordinate 
$(U,(z_1,\cdots,z_n))$ of $M$ with $U\cap\gamma({\Bbb C})\not=\emptyset,\,
Y_i:\gamma^{-1}(U)\to{\Bbb C}$ ($i=1,\cdots,n$) defined by 
$Y_z=\sum\limits_{i=1}^nY_i(z)
(\frac{\partial}{\partial z_i})_{\gamma(z)}$ are holomorphic.  

\vspace{0.5truecm}

\noindent
{\bf Definition 5.} If $Y$ satisfies 
$\nabla^{\Bbb C}_{\gamma_{\ast}(\frac{d}{dz})}
\nabla^{\Bbb C}_{\gamma_{\ast}(\frac{d}{dz})}Y+R^{\Bbb C}(Y,
\gamma_{\ast}(\frac{d}{dz}))\gamma_{\ast}(\frac{d}{dz})=0$, then we call $Y$ 
a {\it complex Jacobi field along} $\gamma$.  
Let $z_0\in{\Bbb C}$.  If there exists a (non-zero) complex Jacobi field $Y$ 
along $\gamma$ with $Y_0=0$ and $Y_{z_0}=0$, then we call $z_0$ a 
{\it complex conjugate radius} of $\gamma(0)$ along $\gamma$.  
Let 
$\delta:{\Bbb C}\times D(\varepsilon)\to M$ be a holomorphic two-parameter map, 
where $D(\varepsilon)$ is the $\varepsilon$-disk centered at $0$ in 
${\Bbb C}$.  Denote by $z$ (resp. $w$) 
the first (resp. second) parameter of $\delta$.  If 
$\delta(\cdot,w_0):{\Bbb C}\to M$ is a complex geodesic for each $w_0\in 
D(\varepsilon)$, then we call $\delta$ a {\it complex geodesic variation}.  

\vspace{0.5truecm}

Easily we can show the following fact.  

\vspace{0.5truecm}

\noindent
{\bf Proposition 8.1.} {\sl Let $\delta:{\Bbb C}\times D(\varepsilon)
\to M$ be a complex geodesic variation.  The complex variational vector field 
$Y:=\delta_{\ast}(\frac{\partial}{\partial w}\vert_{w=0})$ is a complex Jacobi 
field along $\gamma:=\delta(\cdot,0)$.}

\vspace{0.5truecm}

A vector field $X$ on $M$ is said to be {\it real holomorphic} if 
the Lie derivation $L_XJ$ of $J$ with respect to $X$ vanishes.  
It is known that $X$ is a real holomorphic vector field if and only if 
the complex vector field $X-\sqrt{-1}JX$ is holomorphic.  
We have the following fact for a complex Jacobi field.  

\vspace{0.5truecm}

\noindent
{\bf Proposition 8.2.} {\sl Let $\gamma:{\Bbb C}\to M$ be a complex geodesic.  

{\rm(i)} Let $Y$ be a holomorphic vector field along $\gamma$ and $Y_{\Bbb R}$ 
be the real part of $Y$.  Then $Y$ is a complex Jacobi field along $\gamma$ 
if and only if, for any $z_0\in{\Bbb C}$, 
$u\mapsto(Y_{\Bbb R})_{uz_0}$ is a Jacobi field along the geodesic 
$\gamma_{z_0}(\displaystyle{\mathop{\Longleftrightarrow}_{\rm def}}\,\,
\gamma_{z_0}(u):=\gamma(uz_0))$.  

{\rm(ii)} A complex number $z_0$ is a complex conjugate radius of $\gamma(0)$ 
along $\gamma$ if and only if $\gamma(z_0)$ is a conjugate point of 
$\gamma(0)$ along the geodesic $\gamma_{z_0}$.}

\vspace{0.5truecm}

\noindent
{\it Proof.} Let $(z)$ ($z=s+t\sqrt{-1}$) be the natural coordinate of 
${\Bbb C}$.  Let $Y(=Y_{\Bbb R}-\sqrt{-1}JY_{\Bbb R})$ be a holomorphic vector 
field along $\gamma$.  From $L_{Y_{\Bbb R}}J=0$ and $\nabla J=0$, we have 
$$\begin{array}{l}
\hspace{0.6truecm}\displaystyle{\nabla^{\Bbb C}_{\gamma_{\ast}(\frac{d}{dz})}
\nabla^{\Bbb C}_{\gamma_{\ast}(\frac{d}{dz})}Y
+R^{\Bbb C}(Y,\gamma_{\ast}(\frac{d}{dz}))\gamma_{\ast}(\frac{d}{dz})}\\
\displaystyle{=\nabla_{\gamma_{\ast}(\frac{\partial}{\partial s})}
\nabla_{\gamma_{\ast}(\frac{\partial}{\partial s})}Y_{\Bbb R}
+R(Y_{\Bbb R},\gamma_{\ast}(\frac{\partial}{\partial s}))
\gamma_{\ast}(\frac{\partial}{\partial s})}\\
\hspace{0.6truecm}
\displaystyle{-\sqrt{-1}J\left(
\nabla_{\gamma_{\ast}(\frac{\partial}{\partial s})}
\nabla_{\gamma_{\ast}(\frac{\partial}{\partial s})}Y_{\Bbb R}
+R(Y_{\Bbb R},\gamma_{\ast}(\frac{\partial}{\partial s}))
\gamma_{\ast}(\frac{\partial}{\partial s})\right).}
\end{array}\leqno{(8.1)}$$
Assume that $Y$ is a complex Jacobi field.  Then it follows from $(8.1)$ that 
$$\nabla_{\gamma_{\ast}(\frac{\partial}{\partial s})}
\nabla_{\gamma_{\ast}(\frac{\partial}{\partial s})}Y_{\Bbb R}
+R(Y_{\Bbb R},\gamma_{\ast}(\frac{\partial}{\partial s}))
\gamma_{\ast}(\frac{\partial}{\partial s})=0.$$
Let $X:=a\gamma_{\ast}(\frac{\partial}{\partial s})+b\gamma_{\ast}
(\frac{\partial}{\partial t})$ ($a,b\in{\Bbb R}$).  Furthermore, from 
$L_{Y_{\Bbb R}}J=0$ and $\nabla J=0$, we have 
$$\nabla_X\nabla_XY_{\Bbb R}+R(Y_{\Bbb R},X)X=0.$$
Hence we see that 
$u\mapsto(Y_{\Bbb R})_{uz_0}$ is a Jacobi field along 
$\gamma_{z_0}$ for each $z_0\in{\Bbb C}$.  
The converse also is shown in terms of $(8.1), L_{Y_{\Bbb R}}J=0$ and 
$\nabla J=0$ directly.  Thus the statement (i) is shown.  
Assume that $z_0$ is a complex conjugate radius of $\gamma(0)$ along 
$\gamma$.  Then there exists a non-trivial complex Jacobi field $Y$ along 
$\gamma$ with $Y_0=0$ and $Y_{z_0}=0$.  
According to (i), $u\mapsto(Y_{\Bbb R})_{uz_0}$ 
is a Jacobi field along $\gamma_{z_0}$ which vanishes at $u=0,1$.  
Furthermore, it is shown that $u\mapsto(Y_{\Bbb R})_{uz_0}$ is non-trivial 
because so is $Y$.  Hence $\gamma(z_0)$ is a conjugate point of $\gamma(0)$ 
along $\gamma_{z_0}$.  Conversely, assume that $\gamma(z_0)$ is a conjugate 
point of $\gamma(0)$ along $\gamma_{z_0}$.  Then there exists a non-trivial 
Jacobi field $\overline Y$ along $\gamma_{z_0}$ with $\overline Y_0=0$ and 
$\overline Y_1=0$.  There exists the complex Jacobi field $Y$ along $\gamma$ 
with $Y_0=0$ and $\nabla^{\Bbb C}_{{\gamma'}_{z_0}(0)}Y
=\overline Y'_0-\sqrt{-1}J\overline Y'_0$ by the existenceness of solutions 
of a complex ordinary differential equation.  It is easy to show that 
$(Y_{\Bbb R})_{uz_0}=\overline Y_u$ for all $u\in{\Bbb R}$.  Hence we have 
$(Y_{\Bbb R})_{z_0}=\overline Y_1=0$, that is, $Y_{z_0}=0$.  Therefore $z_0$ 
is a complex conjugate radius of $\gamma(0)$ along $\gamma$.  
Thus the statement (ii) is shown.  
\begin{flushright}q.e.d.\end{flushright}

\vspace{0.5truecm}

Next we shall define the notion of the parallel translation along a 
holomorphic curve.  Let $\alpha:D\to(M,J,g)$ be a holomorphic curve, where 
$D$ is an open set of ${\Bbb C}$.  Let $Y$ be a holomorphic vector field along 
$\alpha$.  If $\nabla^{\Bbb C}_{\alpha_{\ast}(\frac{d}{dz})}Y=0$, then we say 
that $Y$ is {\it parallel}.  
For a parallel holomorphic vector field, we can show the following fact.  

\vspace{0.5truecm}

\noindent
{\bf Proposition 8.3.} {\sl Let $\alpha:D\to (M,J,g)$ be a holomorphic curve.  
Take $z_0\in D$ and $v\in(T_{\alpha(z_0)}M)^{(1,0)}$.  Then the following 
statements (i) and (ii) hold.  

(i) There uniquely exists a parallel holomorphic vector field $Y$ along 
$\alpha$ with $Y_{z_0}=v$.  

(ii) Let $Y$ be a holomorphic vector field along $\alpha$ and 
$Y_{\Bbb R}$ be its real part.  Then $Y$ is parallel if and only if, 
for any (real) curve $\sigma$ in $D$, 
$u\mapsto(Y_{\Bbb R})_{\sigma(u)}$ is parallel along $\alpha\circ\sigma$ 
with respect to $\nabla$.}

\vspace{0.5truecm}

\noindent
{\it Proof.} The statement (i) follows from the existenceness and the 
uniqueness of solutions of a complex ordinary differential equation.  The 
statement (ii) is shown as follows.  From $\nabla J=0$ and $L_{Y_{\Bbb R}}J=0$, 
we have 
$\nabla^{\Bbb C}_{\alpha_{\ast}(\frac{d}{dz})}Y=
\frac12(\nabla_{\alpha_{\ast}(\frac{\partial}{\partial s})}Y_{\Bbb R}
-\sqrt{-1}J\nabla_{\alpha_{\ast}(\frac{\partial}{\partial s})}
Y_{\Bbb R})$.  
Hence $Y$ is parallel if and only if $\nabla_{\alpha_{\ast}(\frac{\partial}{\partial s})}Y_{\Bbb R}=0$.  Let $X:=a\gamma_{\ast}(\frac{\partial}{\partial s})+
b\gamma_{\ast}(\frac{\partial}{\partial t})$ ($a,b\in{\Bbb R}$).  
From $\nabla J=0$ and $L_{Y_{\Bbb R}}J=0$, it follows that 
$\nabla_{\alpha_{\ast}(\frac{\partial}{\partial s})}Y_{\Bbb R}=0$ is equivalent 
to $\nabla_XY_{\Bbb R}=0$.  
Therefore, the statement (ii) follows.\hspace{10.4truecm}q.e.d.

\vspace{0.5truecm}

Let $\alpha,\,z_0$ and $v$ be as in the statement of Proposition 8.3.  
There uniquely exists a parallel holomorphic 
vector field $Y$ along $\alpha$ with $Y_{z_0}=v$.  We denote $Y_{z_1}$ by 
$(P_{\alpha})_{z_0,z_1}(v)$.  It is clear that $(P_{\alpha})_{z_0,z_1}$ is 
a ${\Bbb C}$-linear isomorphism of $(T_{\alpha(z_0)}M)^{(1,0)}$ onto 
$(T_{\alpha(z_1)}M)^{(1,0)}$.  We call $(P_{\alpha})_{z_0,z_1}$ 
the {\it parallel translation along} $\alpha$ 
{\it from} $z_0$ {\it to} $z_1$.  

Let $f$ be an immersion of an anti-Kaehler manifold $(M,J,g)$ into another 
anti-Kaehler manifold $(\widetilde M,\widetilde J,\widetilde g)$.  If 
$f_{\ast}\circ J=\widetilde J\circ f_{\ast}$ and $f^{\ast}\widetilde g=g$, 
then we call $f$ an {\it anti-Kaehler immersion} and $(M,J,g)$ an 
{\it anti-Kaehler submanifold immersed by} $f$.  In the sequel, we omit 
the notation $f_{\ast}$.  
In [15], we introduced the notion of a complex focal radius of an anti-Kaehler 
submanifold.  Now we shall define this notion in terms of a complex Jacobi 
field.  Let $v\in T^{\perp}_{p_0}M$ and $\gamma^{\Bbb C}_v
(:D\to \widetilde M)$ be the (maximal) complex geodesic in 
$(\widetilde M,\widetilde J,\widetilde g)$ with 
$(\gamma^{\Bbb C}_v)_{\ast}((\frac{d}{dz})_0)
=\frac12(v-\sqrt{-1}\widetilde Jv)$, where 
$T^{\perp}_{p_0}M$ is the normal space of $M$ at $p_0$ and $D$ is 
a neighborhood of $0$ in ${\Bbb C}$.  

\vspace{0.5truecm}

\noindent
{\bf Definition 6.} If there exists a complex Jacobi field 
$Y$ along $\gamma^{\Bbb C}_v$ with $Y_0(\not=0)\in(T_{p_0}M)^{(1,0)}$ 
and $Y_{z_0}=0$, 
then we call the complex number $z_0$ a {\it complex focal radius of} $M$ 
{\it along} $\gamma^{\Bbb C}_v$.  

\vspace{0.5truecm}

By imitating the proof of (ii) of Proposition 8.2, we can show 
the following fact.  

\vspace{0.5truecm}

\noindent
{\bf Proposition 8.4.} {\sl A complex number $z_0$ is a complex focal radius 
of $M$ along the normal complex geodesic $\gamma^{\Bbb C}_v$ if and only if 
$\gamma^{\Bbb C}_v(z_0)$ is 
a focal point of $M$ along the normal geodesic 
$(\gamma^{\Bbb C}_v)_{z_0}\,
(\displaystyle{\mathop{\Longleftrightarrow}_{\rm def}}\,\,
(\gamma^{\Bbb C}_v)_{z_0}(u):=\gamma^{\Bbb C}_v(uz_0))$, that is, $z_0$ is 
a complex focal radius in the sense of [15].}

\vspace{0.5truecm}

We consider the case where 
$(\widetilde M,\widetilde J,\widetilde g)$ is an anti-Kaehler symmetric space 
$G^{\Bbb C}/K^{\Bbb C}$ and where 
the anti-Kaehler submanifold $M$ is a subset of 
$G^{\Bbb C}/K^{\Bbb C}$ (hence $f$ is the inclusion map).  
For $v\in(T^{\perp}_{b_0K^{\Bbb C}}M)^{\Bbb C}$, we define ${\Bbb C}$-linear 
transformations $\widehat D^{co}_v$ and $\widehat D^{si}_v$ of 
$(T_{b_0K^{\Bbb C}}(G^{\Bbb C}/K^{\Bbb C}))^{\Bbb C}$ by 
$\widehat D^{co}_v:=b_{0\ast}^{\Bbb C}\circ
\cos(\sqrt{-1}{\rm ad}_{\mathfrak g^{\Bbb C}}^{\Bbb C}
((b_{0\ast}^{\Bbb C})^{-1}v))\circ(b_{0\ast}^{\Bbb C})^{-1}$ and 
$\widehat D^{si}_v:=b_{0\ast}^{\Bbb C}\circ
\frac{\sin(\sqrt{-1}{\rm ad}_{\mathfrak g^{\Bbb C}}^{\Bbb C}
((b_{0\ast}^{\Bbb C})^{-1}v))}
{\sqrt{-1}{\rm ad}_{\mathfrak g^{\Bbb C}}^{\Bbb C}((b_{0\ast}^{\Bbb C})^{-1}v)}
\circ(b_{0\ast}^{\Bbb C})^{-1}$, respectively, where 
${\rm ad}_{\mathfrak g^{\Bbb C}}^{\Bbb C}$ is the complexification of the 
adjoint representation ${\rm ad}_{\mathfrak g^{\Bbb C}}$ of 
$\mathfrak g^{\Bbb C}$.  
If, for each $bK^{\Bbb C}\in M$, $b_{\ast}^{-1}(T^{\perp}_{bK^{\Bbb C}}M)\,
(\subset T_{eK^{\Bbb C}}(G^{\Bbb C}/K^{\Bbb C})\subset\mathfrak g^{\Bbb C})$ 
is a Lie triple system (resp. abelian subspace), that is, 
$\exp^{\perp}(T^{\perp}_{bK^{\Bbb C}}M)$ is totally geodesic 
(resp. flat and totally geodesic), then $M$ is said to {\it have section} 
(resp. {\it have flat section}), where 
$\exp^{\perp}$ is the normal exponential map of $M$.  

\vspace{0.5truecm}

\noindent
{\bf Proposition 8.5.} {\sl Let $M$ be an anti-Kaehler submanifold in 
$G^{\Bbb C}/K^{\Bbb C}$ with section and 
$v\in T^{\perp}_{b_0K^{\Bbb C}}M$.  Set $v_{(1,0)}:=\frac12(v-\sqrt{-1}
\widetilde Jv)$.  A complex number $z_0$ is 
a complex focal radius along $\gamma^{\Bbb C}_v$ if and only if 
$$\left.{\rm Ker}\left(\widehat D^{co}_{z_0v_{(1,0)}}
-\widehat D^{si}_{z_0v_{(1,0)}}\circ(A^{\Bbb C})_{z_0v_{(1,0)}}\right)
\right\vert_{(T_{b_0K^{\Bbb C}}M)^{(1,0)}}\not=\{0\},$$
where $A^{\Bbb C}$ is the 
complexification of the shape tensor $A$ of $M$.}

\vspace{0.5truecm}

\noindent
{\it Proof.} Denote by $\widetilde{\nabla}$ (resp. $\widetilde R$) the 
Levi-Civita connection (resp. the curvature tensor) of $G^{\Bbb C}/K^{\Bbb C}$ 
and by $\widetilde{\nabla}^{\Bbb C}$ (resp. $\widetilde R^{\Bbb C}$) their 
complexification.  
Let $Y$ be a holomorphic vector field along $\gamma^{\Bbb C}_v$.  Define 
$\widehat Y:D\to(T_{b_0K^{\Bbb C}}(G^{\Bbb C}/K^{\Bbb C}))^{(1,0)}$ by 
$\widehat Y_z:=(P_{\gamma^{\Bbb C}_v})_{z,0}(Y_z)$ ($z\in D$), where 
$D$ is the domain of $\gamma^{\Bbb C}_v$.  Easily we can show 
$\widetilde{\nabla}^{\Bbb C}_{(\gamma^{\Bbb C}_v)_{\ast}(\frac{d}{dz})}
\widetilde{\nabla}^{\Bbb C}_{(\gamma^{\Bbb C}_v)_{\ast}(\frac{d}{dz})}Y
=(P_{\gamma^{\Bbb C}_v})_{0,z}(\frac{d^2\widehat Y}{dz^2})$.  From 
$\widetilde{\nabla}\widetilde R=0$ (hence $\widetilde{\nabla}^{\Bbb C}
\widetilde R^{\Bbb C}=0$), we have $\widetilde R^{\Bbb C}(Y,
(\gamma^{\Bbb C}_v)_{\ast}(\frac{d}{dz}))(\gamma^{\Bbb C}_v)_{\ast}
(\frac{d}{dz})=(P_{\gamma^{\Bbb C}_v})_{0,z}(R^{\Bbb C}_{b_0K^{\Bbb C}}
(\widehat Y_z,v_{(1,0)})v_{(1,0)})$.  Hence $Y$ is a complex Jacobi field if 
and only if $\frac{d^2\widehat Y}{dz^2}+R^{\Bbb C}_{b_0K^{\Bbb C}}
(\widehat Y_z,v_{(1,0)})v_{(1,0)}=0$ holds.  By noticing 
$$R^{\Bbb C}_{b_0K^{\Bbb C}}(\widehat Y_z,v_{(1,0)})v_{(1,0)}=
-(b_{0\ast}^{\Bbb C}\circ{\rm ad}_{\mathfrak g^{\Bbb C}}^{\Bbb C}
((b_{0\ast}^{\Bbb C})^{-1}v_{(1,0)})^2\circ(b_{0\ast}^{\Bbb C})^{-1})
(\widehat Y_z)$$
and solving this complex ordinary differential equation, we have 
$$\widehat Y_z=\widehat D^{co}_{zv_{(1,0)}}(Y_0)+z\widehat D^{si}_{zv_{(1,0)}}
\left(\left.\frac{d\widehat Y}{dz}\right\vert_{z=0}\right).$$
Since $M$ has section, both 
$\widehat D^{co}_{zv_{(1,0)}}$ and $\widehat D^{si}_{zv_{(1,0)}}$ preserve 
$(T_{b_0K^{\Bbb C}}M)^{\Bbb C}$ (and hence also 
$(T^{\perp}_{b_0K^{\Bbb C}}M)^{\Bbb C}$) invariantly.  
Hence, if $Y_0(\not=0)\in (T_{b_0K^{\Bbb C}}M)^{\Bbb C}$ and $Y_{z_0}=0$ for 
some $z_0$, then we have $\frac{d\widehat Y}{dz}\vert_{z=0}\in(T_{b_0K^{\Bbb C}}
M)^{\Bbb C}$, that is, 
$\frac{d\widehat Y}{dz}\vert_{z=0}=-(A^{\Bbb C})_{v_{(1,0)}}(Y_0)$.  
Hence we have 
$$Y_z=(P_{\gamma^{\Bbb C}_v})_{0,z}((\widehat D^{co}_{zv_{(1,0)}}-
\widehat D^{si}_{zv_{(1,0)}}\circ(A^{\Bbb C})_{zv_{(1,0)}})(Y_0)).
\leqno{(8.2)}$$
From this fact, the statement of this theorem follows.  \hspace{5truecm}q.e.d.

\vspace{0.5truecm}

Let $f:(M,g)\hookrightarrow(\widetilde M,\widetilde g)$ be a 
$C^{\omega}$-isometric immersion between $C^{\omega}$-pseudo-Riemannian 
manifolds and $f^{\Bbb C}:((M_g^{\Bbb C})'_{A,f:i},g_A)\hookrightarrow
(({\widetilde M}_{\widetilde g}^{\Bbb C})_A,\widetilde g_A)$ be 
its complexification (see Definition 2).  

\vspace{0.5truecm}

\noindent
{\bf Definition 7.} For each normal vector $v(\not=0)$ of $M$ 
(in $\widetilde M$), we call a complex focal radius of 
$(M_g^{\Bbb C})'_{A,f:i}$ along 
$\gamma^{\Bbb C}_v$ a {\it complex focal radius of} $M$ 
{\it along the normal geodesic} $\gamma_v$ (in $\widetilde M$).  

\vspace{0.5truecm}

We consider the case where $(\widetilde M,\widetilde g)$ is a 
Riemannian symmetric space $G/K$ of non-compact type and where $M$ has 
section.  
Let $v\in T_{b_0K^{\Bbb C}}^{\perp}M$ and 
$z(=s+t\sqrt{-1})\in{\Bbb C}$.  In [15], we defined the linear map 
$D^{co}_{zv}$ (resp. $D^{si}_{zv}$) of $T_{b_0K^{\Bbb C}}(M_g^{\Bbb C})
(=(T_{b_0K}M)^{\Bbb C})$ into 
$T_{b_0K^{\Bbb C}}(G^{\Bbb C}/K^{\Bbb C})(=(T_{b_0K}(G/K))^{\Bbb C})$ by 
$$\begin{array}{c}
\displaystyle{D^{co}_{zv}:=b_{0\ast}\circ\cos\left(
\sqrt{-1}{\rm ad}_{\mathfrak g^{\Bbb C}}(b_{0\ast}^{-1}
(sv+t\widetilde Jv))\right)\circ b_{0\ast}^{-1}}\\
\displaystyle{\left({\rm resp.}\,\,\,\,
D^{si}_{zv}:=b_{0\ast}\circ\frac{\sin\left(
\sqrt{-1}{\rm ad}_{\mathfrak g^{\Bbb C}}(b_{0\ast}^{-1}
(sv+t\widetilde Jv))\right)}
{\sqrt{-1}{\rm ad}_{\mathfrak g^{\Bbb C}}(b_{0\ast}^{-1}
(sv+t\widetilde Jv))}\circ b_{0\ast}^{-1}\right).}
\end{array}$$
The relations between these operators and the above operators 
$\widehat D^{co}_{zv}$ and $\widehat D^{si}_{zv}$ are as follows:
$$\widehat D^{co}_{zv_{(1,0)}}(X-\sqrt{-1}JX)=D^{co}_{zv}(X)
-\sqrt{-1}J(D^{co}_{zv}(X))
\leqno{(8.3)}$$
and 
$$\widehat D^{si}_{zv_{(1,0)}}(X-\sqrt{-1}JX)=D^{si}_{zv}(X)
-\sqrt{-1}J(D^{si}_{zv}(X)),\leqno{(8.4)}$$
where $X\in T_{b_0K^{\Bbb C}}(M_g^{\Bbb C})$.  
From $(8.2),\,(8.3)$ and $(8.4)$, we have 
$$(Y_{\Bbb R})_z=(P_{(\gamma^{\Bbb C}_v)_z})_{0,1}((D^{co}_{zv}-
D^{si}_{zv}\circ A^{\Bbb C}_{zv})((Y_{\Bbb R})_0)) \leqno{(8.5)}$$
for a complex Jacobi field $Y$ along $\gamma_v^{\Bbb C}$ such that 
$Y_0$ and $\nabla_{(\gamma_v^{\Bbb C})_{\ast}((\frac{d}{dz})_0)}Y$ belong to 
$(T_{b_0K^{\Bbb C}}(M_g^{\Bbb C}))^{\Bbb C}$, where 
$(P_{(\gamma^{\Bbb C}_v)_z})_{0,1}$ is the parallel translation along 
$(\gamma^{\Bbb C}_v)_z$ ($:u\mapsto \gamma_v^{\Bbb C}(uz)$) from $0$ to $1$ and 
$A$ is the shape tensor of $(M,g)$.  Hence we have the following fact.  

\vspace{0.5truecm}

\noindent
{\bf Proposition 8.6.} {\sl Let $M$ be a $C^{\omega}$-Riemannian submanifold 
in a Riemmannian symmetric space $G/K$ of non-compact type.  Then 
$z(\in\Bbb C)$ is a complex focal radius along $\gamma_v$ (in the sense of 
Definition 7) if and only if 
${\rm Ker}(D^{co}_{zv}-D^{si}_{zv}\circ A_{zv}^{\Bbb C})\not=\{0\}$, where 
$A$ is the shape tensor of $M$, that is, $z$ is a complex focal radius along 
$\gamma_v$ in the sense of [14].}

\section{Complex equifocal submanifolds and isoparametric ones} 
In [15], we defined the notion of a complex equifocal submanifold in 
a Riemannian symmetric space of non-compact type by imposing the condition 
related to complex focal radii.  
In the previous section, we defined the notion of a complex focal radius 
for $C^{\omega}$-pseudo-Riemannian submanifold in a general 
$C^{\omega}$-pseudo-Riemannian manifold.  
By imposing the same condition related to complex focal radii, 
we shall define the notion of a complex equifocal submanifold in 
a pseudo-Riemannian homogeneous space.  
Let $M$ be a $C^{\omega}$-pseudo-Riemannian submanifold in a 
$C^{\omega}$-pseudo-Riemannian homogeneous space $\widetilde M$.  
If $M$ has flat section, if the normal holonomy group of $M$ is trivial and 
if, for any parallel normal vector field $v$ of $M$, the complex focal radii 
along $\gamma_{v_x}$ are independent of the choice of $x\in M$ 
(considering their multiplicities), then we call $M$ a {\it complex equifocal 
submanifold}.  
If $M$ has flat section, if the normal holonomy group of $M$ is trivial and 
if, any sufficiently close parallel submanifolds of $M$ have constant mean 
curvature with respect to the radial direction, then $M$ is called an 
{\it isoparametric submanifold with flat section}.  
If, for each normal vector $v$ of $M$, the Jacobi operator $R(\cdot,v)v$ 
preserves $T_xM$ ($x\,:\,$the base point of $v$) invariantly and 
$[A_v,R(\cdot,v)v\vert_{T_xM}]=0$, then $M$ is called a 
{\it curvature-adapted submanifold}, where $R$ is the curvature tensor of 
$\widetilde M$ and $A$ is the shape tensor of $M$.  
By imitating the proof of Theorem 15 in [15], 
we can show the following 
facts for pseudo-Riemannian submanifolds in a semi-simple pseudo-Riemannian 
symmetric space.  

\vspace{0.3truecm}

\noindent
{\bf Proposition 9.1.} {\sl Let $(M,g)$ be a $C^{\omega}$-pseudo-Riemannian 
submanifold in a semi-simple pseudo-Riemannian symmetric space $G/K$ equipped 
with the metric $\widetilde g$ induced from the Killing form of 
$\mathfrak g:={\rm Lie}\,G$.  
Then the following statements {\rm (i)} and {\rm (ii)} hold:

{\rm (i)} If $M$ is an isoparametric submanifold with flat section, then it is 
complex equifocal.  

{\rm (ii)} Let $M$ be a curvature-adapted complex equifocal submanifold.  
If, for any normal vector $w$ of $M$, 
$R^{\Bbb C}(\cdot,w)w\vert_{(T_xM)^{\Bbb C}}$ 
($x:$ the base point of $w$) and the complexified shape operator 
$A^{\Bbb C}_w$ are diagonalizable, then it is an isoparametric submanifold with 
flat section.}

\vspace{0.3truecm}

\noindent
{\it Proof.} Let $M$ be a $C^{\omega}$-pseudo-Riemannian submanifold with 
flat section in $G/K$ whose normal holonomy group is trivial.  Let $v$ be 
a parallel normal vector field on $M$.  Since $M$ has flat section, 
$R(\cdot,v_x)v_x$ preserves $T_xM$ invariantly for for each $x\in M$.  Hence 
the ${\Bbb C}$-linear transformations $D^{co}_{zv_x}$ and $D^{si}_{zv_x}$ 
preserve $(T_xM)^{\Bbb C}(=T_x(M_g^{\Bbb C}))$ invariantly.  Let 
$\eta_{sv}:=\exp^{\perp}\circ sv\,(M\to G/K)$ 
and $M_{sv}:=\eta_{sv}(M)$, where $s$ is sufficiently close to zero.  
Define a function $F_{sv}$ on $M$ by $\eta_{sv}^{\ast}\omega_{sv}=F_{sv}
\omega$, where $\omega$ (resp. $\omega_{sv}$) is the volume element of $M$ 
(resp. $M_{sv}$).  Set $\widehat F_{v_x}(s):=F_{sv}(x)$ ($x\in M$).  From 
$(8.5)$, it follows that $\widehat F_{v_x}$ ($x\in M$) has holomorphic 
extension (which is denoted by $\widehat F_{v_x}^h$) and that 
$$\widehat F^h_{v_x}(z)={\rm det}(D^{co}_{zv_x}-D^{si}_{zv_x}\circ 
A^{\Bbb C}_{zv_x})\quad(z\in{\Bbb C}),\leqno{(9.1)}$$
where $A^{\Bbb C}$ is the complexification of the shape tensor $A$ of $M$, 
that is, the shape tensor of $M_g^{\Bbb C}$ and 
$D^{co}_{zv_x}-D^{si}_{zv_x}\circ A^{\Bbb C}_{zv_x}$ is regarded as a 
${\Bbb C}$-linear transformation of $(T_xM)^{\Bbb C}$.  
By imitating the proof of Corollary 2.6 of [11], 
$M$ is an isoparametric submanifold with flat section if and only if the projection from $M$ to any 
(sufficiently close) parallel submanifold along the sections is volume 
preserving up to a constant factor (i.e., $\widehat F^h_{v_x}$ is independent 
of the choice of $x\in M$ for every parallel normal vector field $v$ of $M$).  
On the other hand, the complex focal 
radii along the geodesic $\gamma_{v_x}$ are catched as zero points of 
$\widehat F^h_{v_x}$.  Hence we see that $M$ is complex equifocal if and only 
if $(F^h_{v_x})^{-1}(0)$ is independent of the choice of $x\in M$ for every 
parallel normal vector field $v$ of $M$.  From these facts, the statement 
(i) follows.  Next we shall show the statement (ii).  
Let $M$ be a curvature-adapted complex equifocal submanifold satisfying the 
conditions of the statement (ii), $v$ be any parallel normal vector field of 
$M$ and $x$ be any point of $M$.  Since $M$ is curvature-adapted, 
$R^{\Bbb C}(\cdot,v_x)v_x$ preserves $(T_xM)^{\Bbb C}$ invariantly, 
$R^{\Bbb C}(\cdot,v_x)v_x\vert_{(T_xM)^{\Bbb C}}$ commutes with 
$A^{\Bbb C}_{v_x}$.  Also, $R^{\Bbb C}(\cdot,v_x)v_x\vert_{(T_xM)^{\Bbb C}}$ and 
$A^{\Bbb C}_{v_x}$ are diagonalizable by the assumption.  Hence they are 
simultaneously diagonalizable.  
Hence, for each $x_0\in M$, there exists a continuous orthonormal tangent 
frame field $(e_1,\cdots,e_n)$ of $(TM)^{\Bbb C}$ defined on a connected open 
neighborhood $U$ of $x_0$ in $M$ such that 
$R^{\Bbb C}(e_i,v)v=-\beta_i^2e_i$ and $A^{\Bbb C}_ve_i=\lambda_ie_i$ 
($i=1,\cdots,n$), where $n:={\rm dim}\,M$, $\beta_i$ and $\lambda_i$ 
($i=1,\cdots,n$) are continuous complex-valued functions on $U$.  
From $(9.1)$, we have 
$$\displaystyle{\widehat F^h_{v_x}(z)=\mathop{\Pi}_{i=1}^n
\left(\cos(\sqrt{-1}z\beta_i(x))-\frac{\lambda_i(x)\sin(\sqrt{-1}z\beta_i(x))}
{\sqrt{-1}\beta_i(x)}\right) \,\,\,\,(x\in U).}
\leqno{(9.2)}$$
Hence we have 
$$\displaystyle{(\widehat F^h_{v_x})^{-1}(0)=\mathop{\cup}_{i=1}^n
\left\{z\,\left\vert\,\cos(\sqrt{-1}z\beta_i(x))=\frac{\lambda_i(x)
\sin(\sqrt{-1}z\beta_i(x))}{\sqrt{-1}\beta_i(x)}\right.\right\}} \leqno{(9.3)}$$
($x\in U$).  Since $M$ is complex equifocal, we have 
$(\widehat F^h_{v_x})^{-1}(0)$ is independent of the choice of $x\in U$.  
Hence, it follows from $(9.3)$ that $\beta_i$ and $\lambda_i$ ($i=1,\cdots,n$) 
are constant on $U$.  
Furthermore, it follows from $(9.2)$ that $\widehat F^h_{v_x}$ is independent 
of the choice of $x\in U$.  From the arbitariness of $x_0$, 
$\widehat F^h_{v_x}$ is independent of the choice of $x\in M$.  
Thus $M$ is an isoparametric submanifold with flat section.  
\hspace{9truecm}q.e.d.

\vspace{0.3truecm}

According to Theorem A of [20], we have the following fact.  

\vspace{0.3truecm}

\noindent
{\bf Proposition 9.2([20]).} {\sl Let $G/K$ be a (semi-simple) 
pseudo-Riemannian symmetric space and $H$ be a symmetric subgroup of $G$, 
$\tau$ (resp. $\sigma$) be an involution of $G$ with $({\rm Fix}\,\tau)_0
\subset K\subset{\rm Fix}\,\tau$ (resp. $({\rm Fix}\,\sigma)_0\subset H
\subset{\rm Fix}\,\sigma$), $L:=({\rm Fix}(\sigma\circ\tau))_0$ and 
$\mathfrak l:={\rm Lie}\,L$, where ${\rm Fix}(\cdot)$ is the fixed point group 
of $(\cdot)$ and ${\rm Fix}(\cdot)_0$ is the identity component of 
${\rm Fix}(\cdot)$.  Assume that $\sigma\circ\tau=\tau\circ\sigma$.  Let $M$ 
be a principal orbit of the $H$-action on $G/K$ through a point $\exp_G(v)K$ 
($v\in \mathfrak q_K\cap\mathfrak q_H$ s.t. ${\rm ad}(v)\vert_{\mathfrak l}
\,:\,$semi-simple), where $\mathfrak q_K:={\rm Ker}(\tau+{\rm id})$ and 
$\mathfrak q_H:={\rm Ker}(\sigma+{\rm id})$.  Then $M$ is a curvature-adapted 
complex equifocal submanifold and, 
for each normal vector $w$ of $M$, $R^{\Bbb C}(\cdot,w)w\vert_{(T_xM)^{\Bbb C}}$ 
($x:$ the base point of $w$) and $A_w^{\Bbb C}$ are diagonalizable.  
Also the orbit $H(eK)$ is a reflective focal submanifold of $M$.}

\vspace{0.3truecm}

By using Theorem 6.1, Propositions 9.1 and 9.2, we prove the following fact.  

\vspace{0.3truecm}

\noindent
{\bf Theorem 9.3.} 
{\sl Let $(G/K,g)$ be a (semi-simple) pseudo-Riemannian symmetric 
space.  Then $((G/K)_g^{\Bbb C})_A$ is invariant with 
respect to the $G$-action on $T(G/K)$ and almost all principal orbits of this 
action are curvature-adapted isoparametric submanifolds with flat section in 
the anti-Kaehler manifold $(((G/K)_g^{\Bbb C})_A,g_A)$ such that the shape 
operators are complex diagonalizable.  Also, the $0$-section($=G/K$) is a 
reflective focal submanifold of such principal orbits.}

\vspace{0.3truecm}

\noindent
{\it Proof.} Since $G$ is a symmetric subgroup of $G^{\Bbb C}$ and 
the involutions associated with $G$ and $K^{\Bbb C}$ commute, it follows from 
Proposition 9.2 that almost all principal orbits of the $G$-action on 
$G^{\Bbb C}/K^{\Bbb C}$ are 
curvature-adapted complex equifocal submanifold such that, for each normal 
vector $w$ of $M$, $R^{\Bbb C}(\cdot,w)w\vert_{(T_xM)^{\Bbb C}}$ 
($x:$ the base point of $w$) and $A_w^{\Bbb C}$ are diagonalizable.  
Also $G(eK^{\Bbb C})(=G/K\subset G^{\Bbb C}/K^{\Bbb C})$ is a reflective focal submanifold of such principal 
orbits.  By Proposition 9.1, such principal orbits are isoparametric 
submanifolds with flat section.  For $g\in G$ and $v\in
((G/K)_g^{\Bbb C})_A\cap T_p(G/K)$, we have 
$$\Phi(g_{\ast}v)=\exp_{g(p)}(\widetilde J_{g(p)}(g_{\ast}v))
=g(\exp_p(\widetilde J_pv))=g(\Phi(v)),$$
where $\Phi$ is as in Section 6 and $\widetilde J$ is the complex structure of 
$G^{\Bbb C}/K^{\Bbb C}$.  Thus $\Phi$ maps the $G$-orbits on 
$((G/K)_g^{\Bbb C})_A$ 
onto the $G$-orbits on $G^{\Bbb C}/K^{\Bbb C}$.  Hence, since 
$\Phi\vert_{((G/K)_g^{\Bbb C})_A}$ is an isometry by Theorem 6.1, almost all 
principal orbits of the $G$-action on $((G/K)_g^{\Bbb C})_A$ are 
curvature-adapted 
isoparametric submanifolds with flat section and their shape operators are 
complex diagonalizable and the $0$-section ($=G/K$) is a reflective focal 
submanifold of such principal orbits.  \hspace{11.15truecm}q.e.d.


\vspace{1truecm}

\noindent
{\Large\bf Concluding remark}

\vspace{0.5truecm}

We shall list up notations used in this paper.  

\vspace{0.5truecm}

$$\begin{tabular}{|c|c|}
\hline
{\scriptsize$J^{\nabla}$} & {\scriptsize the adapted complex structure 
of $\nabla$}\\
\hline
{\scriptsize$J^g$} & {\scriptsize the adapted complex structure of $g$}\\
\hline
{\scriptsize$g_A$} & {\scriptsize the anti-Kaehler metric ass. with $J^g$}\\
\hline
{\scriptsize$M_{\nabla}^{\Bbb C}$} & {\scriptsize the domain of $J^{\nabla}$}\\
\hline
{\scriptsize$M_g^{\Bbb C}$} & {\scriptsize the domain of $J^g$}\\
\hline
{\scriptsize$(M_g^{\Bbb C})_f$} & {\scriptsize the domain of $f^{\Bbb C}$}\\
\hline
{\scriptsize$(M_g^{\Bbb C})_{f:i}$} & {\scriptsize the domain such that 
$f^{\Bbb C}$ is an immersion}\\
\hline
{\scriptsize$(M_g^{\Bbb C})_A$} & {\scriptsize the domain of $(J^g,g_A)$}\\
\hline
{\scriptsize$(M_g^{\Bbb C})_{A,f:i}$} & {\scriptsize $(M_g^{\Bbb C})_A\cap 
(M_g^{\Bbb C})_{f:i}$}\\
\hline
{\scriptsize$(M_g^{\Bbb C})'_{A,f:i}$} & {\scriptsize $(M_g^{\Bbb C})_{A,f:i}
\cap(f^{\Bbb C})^{-1}((\widetilde M_{\widetilde g}^{\Bbb C})_A)$}\\
\hline
\end{tabular}$$


\vspace{0.7truecm}

$$\left(
\begin{array}{c}
{\scriptsize\displaystyle{\nabla\,:\,C^{\omega}-{\rm Koszul}\,\,
{\rm connection}\,\,{\rm of}\,\,M}}\\
{\scriptsize\displaystyle{g\,:\,C^{\omega}-{\rm pseudo-Riemannian}\,\,
{\rm metric}\,\,{\rm of}\,\,M}}\\
{\scriptsize\displaystyle{f\,:\,C^{\omega}-{\rm isometric}\,\,
{\rm immersion}\,\,{\rm of}\,\,(M,g)\,\,{\rm into}\,\,
(\widetilde M,\widetilde g)}}
\end{array}
\right)$$


\vspace{0.7truecm}

\noindent
{\Large\bf References}

\vspace{0.3truecm}

{\small
\noindent
[1] A. Borowiec, M. Ferraris, M. Francaviglia and I. Volovich, 
Almost-complex and almost-product 

Einstein manifolds from a variational principle, 
J Math. Physics {\bf 40} (1999) 3446-3464.

\noindent
[2] A. Borowiec, M. Francaviglia and I. Volovich, 
Anti-K$\ddot a$hlerian manifolds, Differential geom. and 

Its Appl. {\bf 12} (2000) 281-289.

\noindent
[3] M. Cahen and M. Parker, Pseudo-riemannian symmetric spaces, Memoirs 
of the Amer. Math. 

Soc. {\bf 24} No. 229 (1980).  

\noindent
[4] L. Geatti, 
Invariant domains in the complexfication of a noncompact Riemannian symmetric 

space, J. of Algebra {\bf 251} (2002) 619-685.

\noindent
[5] L. Geatti, 
Complex extensions of semisimple symmetric spaces, manuscripta math. {\bf 120} 
(2006) 

1-25.

\noindent
[6] L. Geatti and C. Gorodski, Polar representations of real reductive 
algebraic groups, J of Algebra 

{\bf 320} (2008) 3036-3061.  

\noindent
[7] J. Hahn, Isoparametric hypersurfaces in the pseudo-Riemannian 
space forms, Math. Z. {\bf 187} 

(1984) 195-208.  

\noindent
[8] J. Hahn, Isotropy representations of semisimple symmetric spaces 
and homogeneous hypersur-

faces, J. Math. Soc. Japan {\bf 40} (1988) 271-288.  

\noindent
[9] B. C. Hall and W. D. Kirwin, Adapted complex structures and the 
geodesic flow, Math. Ann. 

{\bf 350} (2011) 455-474.

\noindent
[10] S. Helgason, 
Differential geometry, Lie groups and symmetric spaces, Pure Appl. Math. 80, 

Academic Press, New York, 1978.

\noindent
[11] E. Heintze, X. Liu and C. Olmos, Isoparametric submanifolds and 
a Chevalley type restricti-

on theorem, Integrable systems, geometry, and topology, 151-190, 
AMS/IP Stud. Adv. Math. 

36, Amer. Math. Soc., Providence, RI, 2006.

\noindent
[12] S. Kobayashi and K. Nomizu, 
Foundations of differential geometry, Interscience Tracts in 

Pure and Applied Mathematics 15, Vol. I, New York, 1969.

\noindent
[13] S. Kobayashi and K. Nomizu, 
Foundations of differential geometry, Interscience Tracts in Pure 

and Applied Mathematics 15, Vol. II, New York, 1969.

\noindent
[14] N. Koike, 
Submanifold geometries in a symmetric space of non-compact 
type and a pseudo-

Hilbert space, Kyushu J. Math. {\bf 58} (2004), 167--202.

\noindent
[15] N. Koike, 
Complex equifocal submanifolds and infinite dimensional anti-
Kaehlerian isopara-

metric submanifolds, Tokyo J. Math. {\bf 28} (2005), 
201--247.

\noindent
[16] N. Koike, 
Actions of Hermann type and proper complex equifocal submanifolds, 
Osaka J. 

Math. {\bf 42} (2005) 599-611.

\noindent
[17] N. Koike, 
A splitting theorem for proper complex equifocal submanifolds, Tohoku Math. J. 

{\bf 58} (2006) 393-417.

\noindent
[18] N. Koike, A Chevalley type restriction theorem for a proper complex 
equifocal submanifold, 

Kodai Math. J. {\bf 30} (2007) 280-296.

\noindent
[19] N. Koike, The homogeneous slice theorem for the complete complexification 
of a proper 

complex equifocal submanifold, Tokyo J. Math. 
{\bf 33} (2010) 1-29.

\noindent
[20] N. Koike, Hermann type actions on a pseudo-Riemannian symmetric space, 
Tsukuba J. Math. 

{\bf 34} (2010) 137-172.

\noindent
[21] L. Lempert and R. Sz$\ddot o$ke, Global solutions of the homogeneous 
complex Monge-Amp$\grave e$re eq-

uation and complex structures on the tangent bundle of 
Riemannian manifolds, Math. Ann. 

{\bf 290} (1991) 689-712.

\noindent
[22] B. O'Neill, 
Semi-Riemannian Geometry, with Applications to Relativity, 
Academic Press, 

New York, 1983.

\noindent
[23] T. Ohshima and J. Sekiguchi, The restricted root system of 
a semi-simple symmetric pair, 

Advanced Studies in Pure Mathematics {\bf 4} (1984) 433-497.

\noindent
[24] D. C. Robinson, The real geometry of holomorphic four-metrics, 
J. Math. Physics {\bf 43} (2002) 

2015-2028.

\noindent
[25] W. Rossmann, 
The structure of semisimple symmetric spaces, Can. J. Math. {\bf 1} 
(1979), 157

--180.

\noindent
[26] R. Sz$\ddot{{{\rm o}}}$ke, Complex structures on tangent 
bundles of Riemannian manifolds, Math. Ann. {\bf 291} 

(1991) 409-428.

\noindent
[27] R. Sz$\ddot{{{\rm o}}}$ke, Automorphisms of certain Stein 
manifolds, Math. Z. {\bf 219} (1995) 357-385.

\noindent
[28] R. Sz$\ddot{{{\rm o}}}$ke, Adapted complex structures and 
geometric quantization, Nagoya Math. J. {\bf 154} 

(1999) 171-183.

\noindent
[29] R. Sz$\ddot{{{\rm o}}}$ke, Involutive structures on the 
tangent bundle of symmetric spaces, 
Math. Ann. {\bf 319} 

(2001), 319--348.

\noindent
[30] R. Sz$\ddot{{{\rm o}}}$ke, Canonical complex structures associated to 
connections and complexifications of 

Lie groups, Math. Ann. {\bf 329} (2004), 553--591.

\noindent
[31] C. L. Terng and G. Thorbergsson, 
Submanifold geometry in symmetric spaces, J. Differential 

Geom. {\bf 42} (1995), 665--718.

}

\vspace{1truecm}

\begin{flushleft}
{\small\textit{Department of Mathematics, Faculty of Science, 
Tokyo University of Science,}}

{\small\textit{1-3 Kagurazaka Shinjuku-ku, Tokyo 162-8601 Japan}}

{\small\textit{E-mail address}: koike@ma.kagu.tus.ac.jp}
\end{flushleft}

\end{document}